%% file: article.tex
\documentclass[11pt]{article}
\usepackage{epsfig}
\usepackage[utf8]{inputenc}
\usepackage[T1]{fontenc}

\usepackage{color}
\usepackage{graphicx}

\usepackage{amsmath} 
\usepackage{amssymb}

\newcommand{\RR}{\ensuremath{\mathbb R}}
\newcommand{\WW}{\ensuremath{\mathbb W}}

\newcommand{\qp}{\dot{\bf{q}}}
\newcommand{\up}{\tilde{{{\bf u}}}}
\newcommand{\un}{\tilde{{{u}}}_n}
\newcommand{\vp}{{\bf v}}

\newcommand{\no}{{\bf n}}
\newcommand{\st}{\textbf r}
\newcommand{\sn}{r_n}
\newcommand{\ta}{{\mathbf \tau}}
\newcommand{\Db}{{\bf D}^{c,i}}
\newcommand{\Dn}{{D}^{c,i}_n}

\newcommand{\Ab}{{\bf A}^{c,i}}

\newcommand{\At}{{\bf A}^{c,i}_t}

\newcommand{\Zt}{{\bf Z}}

\newtheorem{rmk}{Remark}

\title{On enhanced descend algorithms for solving frictional multi-contact problems : applications to the Discrete Element Method}
\date{\today}
\author{Serge Dumont$^1$}

\begin{document}

\maketitle

{\small
\begin{center}
${^{1}}$ LAMFA, Universit\'e de Picardie Jules Verne -  CNRS UMR 7352, 33, rue Saint-Leu,\\
 80 000 Amiens, France. \\
e-mail : serge.dumont\symbol{'100}u-picardie.fr\\
\end{center}
}

{\bf Abstract}
In this article, we present various numerical methods to solve  multi-contact problems within the Non-Smooth Discrete Element Method.
The techniques considered to solve the frictional unilateral conditions are based both on the bi-potential theory introduced by de Saxc\'e et al. \cite{dSF91}
and the Augmented Lagrangian theory introduced by Alart et al. \cite{AC91}.
Following the ideas of Feng et al. \cite{FJCM05}, a new Newton method is developed to improve these classical algorithms and 
numerical experiments are presented to show that these
methods are faster than the previous ones and provides results with a better quality.
\bigskip


\section{Introduction}
This is a first draft of a paper that will be submitted in a near future.
\medskip

The paper is organized as follow: in the next part, we present the equations to be solved for the Discrete Element Method, and the frictional
contact law considered. In the third part, we first present two classical methods to numerically solve the full problem, the first one based on the bi-potential
theory, and the second one on the Augmented Lagrangian theory. Then, we show how these methods can be enhanced using an appropriate
Newton method. The last part on this article is devoted to the numerical experiments in order to show the main properties of these algorithms.

\section{Problem Setting}




\subsection{The equations of motion of a multi-contact system}
Classically (see for example \cite{J99,JM92,M88}), the motion of a multi-contact system is described using a global 
generalized coordinate $\bf q$ (for $N_p$ particles, ${\bf q}\in \RR^{\tilde d \times N_p}$, where $\tilde d=6$ for a 3D problem and
$\tilde d=3$ for a 2D problem). Due to the possible shocks between particles, the equations of motion has to be formulated
in term of differential measure equation:
\begin{equation}\label{eqmotion}
{\mathbb M} d\qp + {\bf F}^{int}(t,{\bf q},\qp)dt = {\bf F}^{ext}(t,{\bf q},\qp)dt + d{R}
\end{equation}

where
\begin{itemize}
\item $\mathbb M$ represents the generalized mass matrix;

\item $ {\bf F}^{int}$ and $ {\bf F}^{ext}$ represent the internal and external forces respectively;

\item $d{\bf R}$ is a non-negative real measure, representing the reaction forces and impulses between particles in contact.

\end{itemize}

For the sake of simplicity and without lost of generality, only the external forces are considered in
the following. The internal forces are neglected because the general case can be easily derived through a linearizing procedure.
\medskip

Then, for the numerics, the equation (\ref{eqmotion}) is integrated on each time interval $[t_k,t_{k+1}]$, and approximated using a
$\theta$-method with $\theta\in]\frac12,1]$ for stability reason (see \cite{RA04}).

Therefore, the classical approximation of equation (\ref{eqmotion}) yields
\begin{equation}\label{eqdiscr1}
\left\{\begin{array}{l}
\mathbb{M}(\qp_{k+1}-\qp_n) = \Delta t(\theta {\bf F}_{k+1}+(1-\theta){\bf F}_k)+{\bf R}_{k+1} \\
{\bf q}_{k+1}={\bf q}_k+\Delta t \theta \qp_{k+1}+\Delta t(1-\theta)\qp_k
\end{array}
\right.
\end{equation}

We will denote $\qp_k^{free}=\qp_k + \mathbb{M}^{-1}\Delta t(\theta {\bf F}_{k+1}+(1-\theta){\bf F}_k)$
the free velocity (velocity when the contact forces vanish). Then, the first equation in (\ref{eqdiscr1}) becomes
\begin{equation}
\qp_{k+1}=\qp_k^{free}+ \mathbb{M}^{-1}{\bf R}_{k+1}.
\end{equation}

In order to write the contact law, for a contact $c$ between two particles ($1\leq c\leq N_c$, where $N_c$ is
is the total number of contact), we define the local-global mapping 

\begin{equation}
\left\{\begin{array}{l}
{\bf u}^c=P^*({\bf q},c)\qp \\
{\bf R}=P({\bf q},c)\st^c
\end{array}\right.
\end{equation}

where ${\bf u}^c$ is the local relative velocity between the two bodies in contact and $\st^c$ is the local
contact forces (${\bf u}^c,\st^c\in \RR^d$ where $d$ is the dimension of the problem, and $P^*$ is the transpose of matrix $P$).
We also denote ${\mathbb P}({\bf q})$ the total-global mapping, for ${\bf u}$ and $\st$ in
$\RR^{d\times N_c}$ (vectors composed of all relative velocity and contact forces respectively):
\begin{equation}\label{globalmap}
\left\{\begin{array}{l}
{\bf u}={\mathbb P}^*({\bf q})\qp \\
{\bf R}={\mathbb P}({\bf q})\st
\end{array}\right.
\end{equation}

In the discretization, a prediction of $\bf q$ 
is computed to estimate the mapping ${\mathbb P}({\bf q})$ (see equations (\ref{globalmapapp1}) and (\ref{globalmapapp2}) in the following).
\bigskip

Using the equations (\ref{eqdiscr1}) and (\ref{globalmap}), the discretization of the motion of
a multi-contact system, with frictional contact between particles can be written:

\begin{equation}\label{eqdiscr2}
\left\{\begin{array}{l}
\up_{k+1}= \up_k^{free} +{\mathbb W}\st_{k+1} \\
law_c(\up^c_{k+1},\st^c_{k+1})=\mbox{.true.}\qquad \forall c\in \{1,2,...,N_c\}
\end{array}
\right.
\end{equation}

where ${\mathbb W}={\mathbb P}^* {\mathbb M}^{-1} {\mathbb P}$ is the Delassus operator, 
and $\up_k^{free}={\mathbb P}^*\qp_k^{free}$ is the relative free velocity.
Notice that a Newton impact law is also considered (see \cite{M88} and equation (\ref{newtonlaw}) in the following),
that modify ${\bf u}_k$ and ${\bf u}^{free}_k$ by $\up_k$ and $\up^{free}_k$ respectively.
\medskip

The second equation in (\ref{eqdiscr2}) is the implicit frictional contact law that is in our case
the classical Signorini condition and Coulomb's friction law.

\subsection{The frictional contact law}
\begin{figure}[htbp]
\begin{center}
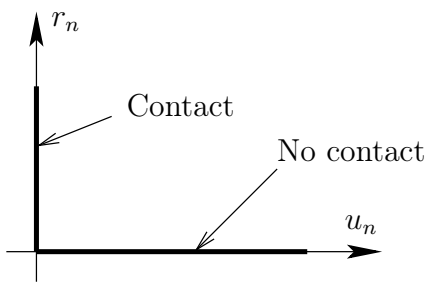
\caption{The Signorini conditions}
\label{signorini}
\end{center}
\end{figure}

In the local coordinates system defined by the local normal vector $\no$ and the tangential vector ${\bf t}\perp\no$, any element
${\bf u}$ and $\st$ can be uniquely decomposed as ${\bf u}=u_n\no+{\bf u}_t$ and $\st=\sn\no+\st_t$ respectively. In these coordinates,
the unilateral contact law can be stated using the Signorini's conditions (see figure \ref{signorini} for a graphical representation):
\begin{equation}
u_n\geq 0,\quad r_n\geq 0,\quad u_nr_n=0.
\end{equation}

\begin{figure}[htbp]
\begin{center}
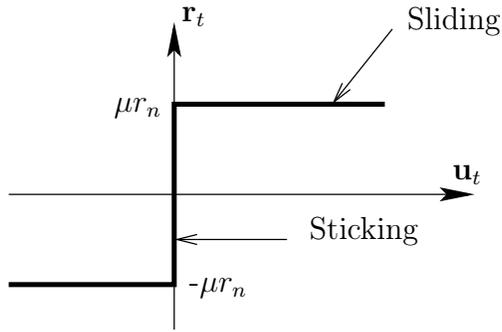
\caption{The Coulomb conditions}
\label{coulomb}
\end{center}
\end{figure}

On the other hand, the Coulomb's law of friction can be stated using the algorithmic form (see figure \ref{coulomb} for a graphical representation):
\begin{equation}
\left[\begin{array}{lll}
\mbox{If } \sn=0\quad &\mbox{ then } u_n\geq 0 & \mbox{! No contact}\\
 & & \\
\mbox{Else if } \sn>0 \mbox{ and } \|\st_t\|<\mu \sn & \mbox{ then } {\bf u}=0 & \mbox{! Sticking}\\
 & & \\
\mbox{Else } \sn>0 \mbox{ and } \|\st_t\|=\mu \sn & \mbox{ then } \exists \lambda \geq 0 \mbox{ such that } {\bf u}_t=\lambda \frac{\st_t}{\|\st_t\|} & \mbox{! Sliding}
\end{array}\right.
\end{equation}

\begin{figure}[htbp]
\begin{center}
\includegraphics[width=6.0cm]{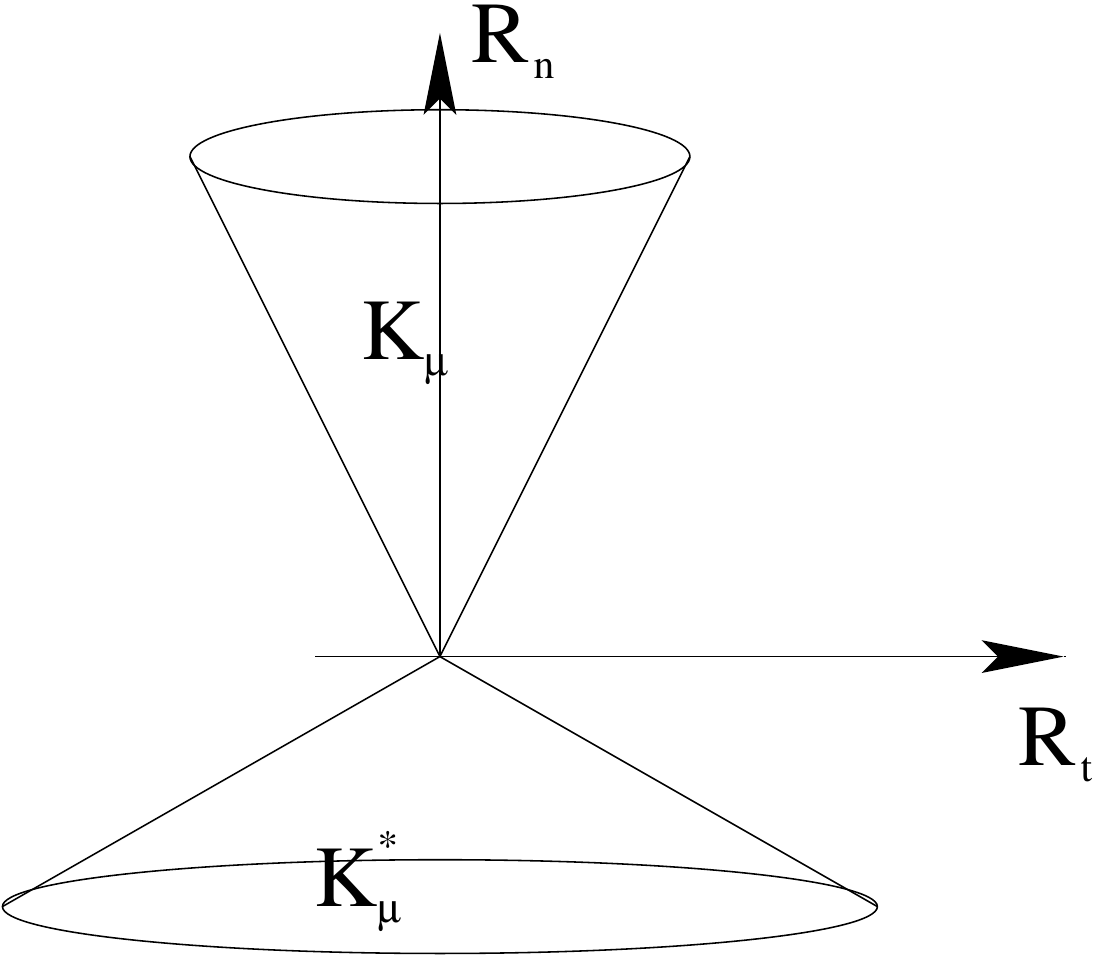}
\caption{The Coulomb's cone}
\label{coulombcone}
\end{center}
\end{figure}

For a given friction coefficient $\mu$, let $K_\mu$ be the isotropic Coulomb's cone, which defines the set of admissible forces (see figure \ref{coulombcone}):
\begin{equation}
K_\mu=\left\{ \st=\sn\no+\st_t :\  \|\st_t\|-\mu \sn\leq 0 \right\}
\end{equation}

The previous law can be also written:
\begin{equation}
\left[\begin{array}{lll}
\mbox{If } \sn=0\quad &\mbox{ then } u_n\geq 0 & \mbox{! No contact}\\
 & & \\
\mbox{Else if } \st\in I(K_\mu) & \mbox{ then } {\bf u}=0 & \mbox{! Sticking}\\
 & & \\
\mbox{Else } \sn>0 \mbox{ and } \st \in B(K_\mu) & \mbox{ then } \exists \lambda \geq 0 \mbox{ such that } {\bf u}_t=\lambda \frac{\st_t}{\|\st_t\|} & \mbox{! Sliding}
\end{array}\right.
\end{equation}

where $I(K_\mu)$ and $B(K_\mu)$ are respectively the interior and the boundary of the cone $K_\mu$.

\section{Numerical Resolution of the contact/friction problems}

We will describe in this section the numerical algorithms that will be considered in the following. Generally, to solve the problem
(\ref{eqdiscr2}), the numerical algorithms considered are based on two levels: the global level where the equations
of motion are solved, and the local level devoted to the resolution of the contact law. 

\subsection{Resolution of the global problem : the Non Linear Gauss Seidel Method (NLGS) }

In this paragraph, we describe the algorithm used at the global level to solve the problem (\ref{eqdiscr2}).
Following the ideas of Jean and Moreau \cite{J99, M88}, we use the non-linear Gauss-Seidel algorithm which is
the most commonly used. It consists in considering successively each contact until the convergence.
The numerical criterion used to state the convergence will be studied latter in the paper.

This method is intrinsically sequential but it is possible to used a simple multi-threading technique which consists
in splitting the contact loop into several threads. This method has been studied in \cite{RDA04} in the case
where the local algorithm is based on the Augmented Lagrangian method.

Notice that it is also possible to consider at this stage more sophisticated method such as a conjugate gradient type method
(see for example \cite{RA04}).


\subsection{The standard bi-potential based method (SBP) }

In this paragraph, we provide a first method to solve the contact problem, at the local level (contact point between two particles).
The method is based on the notion of bi-potential, introduced by de Saxc\'e et al. \cite{dSF91}.

\bigskip

Using the bi-potential framework, it can be shown (see for example \cite{dSF91,F00,Fortin2005,S06}) that a couple
$({\bf u},\st)$ verifies the Signorini-Coulomb contact rules if
\begin{equation}\label{SCrules}
b_c({\bf v},{\bf s})+{\bf v}\cdot {\bf s}\geq b_c({\bf u},\st)+{\bf u}\cdot \st=0\qquad \forall {\bf v},{\bf s}
\end{equation}

where $b_c$ is the bi-potential
\begin{equation}
b_c({-\bf u},\st)=\Psi_{\RR^+}(u_n) +\Psi_{K_\mu}(\st) +\mu \sn\|{\bf }u_t\|
\end{equation}

and $\Psi_C$ stands for the indicatrix function of the set $C$: $\Psi_C(x)=0$ if $x\in C$, $\Psi_C(x)=+\infty$ if $x\notin C$.

Consequently, the contact law can be written in a compact form of an implicit subnormality rule (or a differential inclusion rule):
\begin{equation}
-{\bf u} \in \partial_\st b_c(-{\bf u},\st).
\end{equation}

Then, for a contact $c$, at a NLGS iteration $i$, knowing the relative velocity $\up^{c,i}$, the algorithm to compute  $\st^{c,i+1}$ from $\st^{c,i}$ is based on
the minimization of the bi-potential (see for exemple \cite{F00}, page 51), using the inequality:

\begin{equation}\label{minbipo}
b_c(-\up^{c,i},\st)+\up^{c,i}\cdot\st \geq b_c(-\up^{c,i},\st^{c,i+1})+\up^{k,i}\cdot\st^{c,i+1}\qquad \forall \st\in K_\mu
\end{equation}

or $g(\st)\geq g(\st^{c,i+1})$, $\forall \st\in K_\mu$, if we denote
\begin{equation}\label{defg}
g(\st)=\Psi_{\RR^+}(\un^{c,i})+\Psi_{K_\mu}(\st)+\mu \sn\|\up^{c,i}_t \|+\up^{c,i}\cdot\st.
\end{equation}

The minimization of  (\ref{minbipo}) is classically realized using a projected gradient projection
(Uzawa method) without considering the singular term $\Psi_{\RR^+}(\un^{c,i})$. This minimization can
also be viewed as the proximal point of the augmented force $\st-\rho \up$, with respect to the function 
$\st\mapsto \rho b_c(-\up,\st)$ (see for example \cite{dSF91,F00,Fortin2005}):
\[\st=prox(\st-\rho\up,\rho b_c(-\up,\st)).\]
\medskip

More precisely, the Uzawa method leads to compute the augmented force
$\ta^{c,i+1}=\st^{c,i}-\rho \nabla \tilde{g} (\st^{c,i})$, where $\tilde{g}$ is the differential part of $g$:
$$
\nabla \tilde{g}(\st^{c,i})=\nabla_\st (\mu \sn\|\up^{c,i}_t \|+\up^{c,i}\cdot\st)=\mu \|\up^{c,i}_t \| \no + \up^{c,i},
$$
and to consider the force at next step as a projection of the augmented force onto the set of admissible force $\st^{c,i+1}=proj(\ta^{c,i+1},K_\mu)$, that provides equations (\ref{uzawa1}) and (\ref{uzawa2}) in the resolution algorithm 
of the global problem. The $proj(\ta^{c,i+1},K_\mu)$ stands for the orthogonal projection over the convex $K_\mu$, that can be computed
exactly (see \cite{F00}).
\bigskip

This algorithm will be referred as the SBP (Standard Bi-Potential) method above and throughout.
\bigskip


For a sake of simplicity, we denote hereafter the descent direction
$$\Db=\mu \|\up^{c,i}_t \| \no + \up^{c,i}.$$

\begin{rmk}
A first improvement of this method could be to compute the optimal step $\rho^{c,i}$.
To do so, we have to minimize

\begin{equation}\label{minrho}
\rho\mapsto g(\st^{c,i}-\rho \Db),
\end{equation}

or, more precisely,
\begin{equation}\label{minrho1}
\begin{array}{ll}
\rho\mapsto & \ \ \ \Psi_{\RR^+}(\un^{c,i})+\Psi_{K_\mu}(\st^{c,i}-\rho \Db)+\mu (\sn^{c,i}-\rho \Db\cdot \no) \|\up^{c,i}_t \|+\up^{c,i}\cdot(\st^{c,i}-\rho \Db)\\
 & =\Psi_{\RR^+}(\un^{c,i})+\Psi_{K_\mu}(\st^{c,i}-\rho \Db)-\rho \Db\cdot(\mu \|\up^{c,i}_t \| \no + \up^{c,i})+Cte\\
 & = \Psi_{\RR^+}(\un^{c,i})+\Psi_{K_\mu}(\st^{c,i}-\rho \Db)-\rho \|\Db\|^2+Cte.\\
 \end{array}
\end{equation}
We can observe that this method do not permit to choose an optimal parameter $\rho$ since $g$, as a function of $\rho$,  is linear, excepted in the
case where $ \Db\notin K_\mu$. A solution could be to modify the function $g$, for example by replacing $\up^{c,i}$
by a prediction of $\up^{c,i+1}$ using the equations of the dynamics. Unfortunately, this method do not provides
good numerical results.
\end{rmk}
\medskip\goodbreak

Then, the standard bi-potential based algorithm (SBP) can be written (see \cite{S06} for example):

\begin{itemize}
\item Loop on the step time $k$

\begin{itemize}
\item Prediction of a position (for the computation of the local-global mapping): \\
\begin{equation}\label{globalmapapp1}
{\bf q}_{k+\frac12}={\bf q}_k+\frac{\Delta t}{2}\qp_k;
\end{equation}
\item Initialization of the motion: $\qp^{0}_{k+1}=\qp_{k}^{free}$
(initialization of the contact forces with ${\bf R}=0$).

\item  Loop on $i\geq 0$ (NLGS), until convergence

\begin{itemize}
\item Loop on the contacts $c$:

\begin{itemize}

\item Computation of the local-global mapping
\begin{equation}\label{globalmapapp2}
\dot{\bf  u}^-=P^*({\bf q}_{k+\frac12},c)\qp_k\ ;\qquad\dot{\bf u }^{c,+i}=P^t({\bf q}_{k+\frac12},c)\qp^{i}_{k+1}
\end{equation}
\item Newton shock law
\begin{equation}\label{newtonlaw}
\un^{c,i}=\frac{u^{c,+i}_n+e_nu^-_n}{1+e_n}\ ;\qquad \up^{c,i}_t=\frac{{\bf u}^{c,+i}_t+e_n{\bf u}^-_t}{1+e_t} 
\end{equation}

\item Prediction of the reaction:
\begin{equation}\label{uzawa1}
\ta^{c,i+1}= \st^{c,i} -\rho \left[ \up^{c,i}_t+( \un^{c,i}+\mu \|\up^{c,i}_t \|)\no  \right]
\end{equation}
\item Correction of the reaction:
\begin{equation}\label{uzawa2}
\st^{c,i+1}=proj(\ta^{c,i+1},K_{\mu})
\end{equation}

\item Actualization of the generalized displacement:
\begin{equation}\label{uzawa3}
\displaystyle \qp^{i+1}_{k+1}=\qp_{k}^{free}+{\mathbb M}^{-1}(\sum_{\alpha\leq c}P({\bf q}_{k+\frac12},\alpha)\st^{\alpha,i+1}
 +\sum_{ \alpha>c}P({\bf q}_{k+\frac12},\alpha)\st^{\alpha,i})
\end{equation}
\end{itemize}
\item End of the loop on contacts $c$.\\
\end{itemize}

\item End of the loop on  $i$ of NLGS when the convergence is reached: $\qp_{k+1}=\qp_{k+1}^{i+1}$

\item Actualization of the generalized displacements: ${\bf q}_{k+1}={\bf q}_{k+\frac12}+\frac{\Delta t}{2}\qp_{k+1}$
\end{itemize}
\item End of the loop on the step time $k$.
\end{itemize}

\begin{rmk}
Notice that only one iteration of the Uzawa algorithm at the local level is considered. Various previous studies
(see for example \cite{JF08}) show that there is no significant improvement of the method if several iterations of the Uzawa algorithm
are considered at this stage.
\end{rmk}

\subsection{Newton method and enhanced bi-potential method (EBP)}

We introduce in this section a Newton method in order to speed up the convergence of the computation of the solution.
This method has been already used, especially in the case of the augmented lagrangian method developed
by Alart et al. \cite{AC91}, and the ideas presented in this article follows those of Feng et al. \cite{FJCM05} and 
 have been adapted to the problem of the discrete element method.
The main idea of this technique is to find the solution of the optimization problem, not as a minimum of a functional, 
but rather as a zero of a function, using the Euler equation of the problem. Then a standard Newton method can be developed
to solve this Euler equation.

The technique is first described in the case of the bi-potential framework, and will adapted to the augmented lagrangian method
farther.

\bigskip

We recall that the local problem that has to be solved, for each contact $c$ can be written

\begin{equation}\label{locprob1}
\left\{\begin{array}{l}
\displaystyle \up^{c}_{k+1}=\up^{c,free}_k+\sum_{\alpha=1}^{N_c} W_{c\alpha} \st^\alpha
\ \\
\st^c=proj(\ta^c,K_\mu)
\end{array}\qquad \forall c=1,\ ...,\ N_c
\right.\end{equation}

where $\ta^c=\st^c-\rho(\mu\|\up_t^c\|\no+\up)$ is the augmented reaction (see \ref{uzawa1}), and
$W_{c\alpha}=P^*({\bf q}_{k+\frac12,c}){\mathbb M}^{-1}P({\bf q}_{k+\frac12,\alpha})$ is the local Delassus operator.
\medskip

This problem can be written equivalently
\begin{equation}\label{locprob2}
\left\{\begin{array}{l}
\displaystyle \up^c_{k+1}-\up^{c,free}_{k}-\sum_{\alpha=1}^{N_c} W_{c\alpha} \st^\alpha=0 \\
\ \\
\st^c-proj(\ta^c,K_\mu)=0
\end{array}\right.\qquad \forall c=1,\ ...,\ N_c
\end{equation}

Reminding now that we want to use a Newton algorithm to solve theses equations inside the Non Linear Gauss Seidel loop on the variable $i$, we define now,
for each contact $c=1,\ ...,N_c$, the function
\[
 f^i_c(\chi)=
\left(\begin{array}{c}
\displaystyle\up^{c,i}-\up^{c,free}_{k}-\sum_{\alpha=1}^{N_c}  W_{c\alpha} \st^{\alpha,i} \\
\ \\
\Zt^{c,i} 
\end{array}\right)
\]

where :
\begin{itemize}
\item the vector $\Zt^c$ is the error on the prediction of the reaction
\begin{equation}\label{errorpred}
\Zt^{c,i}(\st^{c,i},\up^{c,i})=\st^{c,i}-proj(\ta^{c,i},K_\mu),
\end{equation} 
\item $\chi_c=(\st^{c,i},\up^{c,i})^t$, 
\item $\chi=(\chi_1,\chi_2,...,\chi_{N_c})^t$
\end{itemize}

\begin{rmk}
The first equality in the relation $f(\chi)=0$ is the equation of motion for the bodies in contact, and the second
relation is the frictional Coulomb law between the bodies in contact, written within the bipotential framework.
\end{rmk}

Then we have to write a Newton algorithm to solve the problem $f(\chi)=0$.
This algorithm can be written, for a contact $c$, by substituting equations (\ref{uzawa1}) and (\ref{uzawa2}) in algorithm (SBP)
by the followings:

\begin{itemize}
\item Initialization:
$$
\chi^0_c=\left( \st^0=\st^{c,i},\ \vp^0=\up^{c,i} \right)^t,\quad \ell=0
$$

\item Loop on $\ell$, until convergence:
\begin{itemize}
\item $\ta^c_\ell=\st^\ell-\rho(\mu\|\vp^\ell_t\|\no+\vp^\ell)$

\item Resolution:
\begin{equation}\label{eqnewt}
\left[ \frac{\partial f_c}{\partial \chi^c}(\chi^\ell)\right] \Delta \chi_c= -f_c(\chi^\ell)
\end{equation}

\item Actualization: $\chi_c^{\ell+1}=\chi_c^\ell+\Delta \chi_c$
\end{itemize}
\item End of the loop on $\ell$ until convergence, $\up^{c,i+1}=\vp^{\ell}$ and $\st^{c,i+1}=\st^\ell$.

\end{itemize}

\begin{rmk}
This algorithm needs more than one iteration at each Non Linear Gauss Seidel iteration to be efficient. As a consequence and compared
to the Uzawa algorithm, the solution in the Newton algorithm is controlled by both the local (iteration $\ell$) and global convergence criteria
(iteration $i$, see \cite{FJCM05,JF08}).
\medskip

The local convergence criterion for the Newton algorithm is defined by:
\begin{equation}\label{er_newt}
\varepsilon^c_{Newt}(\chi_\ell)=\|\vp^\ell-{\bf u}^{c,free}_k-W\st^\ell\|+\|\st^\ell-proj(\st^\ell,K_\mu)\|
\end{equation}
This criterion measure $f_c(\chi_\ell)$ that has to be sufficiently small.
\end{rmk}

The matrix $\left[ \frac{\partial f_c}{\partial \chi_c}(\chi)\right]$ represents the tangential matrix
of the local equations for the contact $c$. This matrix is of dimension  $6\times6$ for a 3 dimensional problem,
and $3\times 3$ for a 2 dimensional problem. For a 3 dimensional problem, the general form of this matrix is the following: 
\begin{equation}
\left[  \frac{\partial f_c}{\partial \chi_c}(\chi)\right]
=
\left[\begin{array}{cc}
-W & Id_{3\times3} \\
A_c & B_c\\
\end{array}\right]
\end{equation}
where
\begin{equation}
A_c=\left[ \frac{\partial Z_c}{\partial \sn}\left|  \frac{\partial Z_c}{\partial r_{t_1}}\right| \frac{\partial Z_c}{\partial r_{t_2}}\right]
\qquad
B_c=\left[ \frac{\partial Z_c}{\partial v_n}\left|  \frac{\partial Z_c}{\partial v_{t_1}}\right| \frac{\partial Z_c}{\partial v_{t_2}}\right]
\end{equation}
\goodbreak

The matrices $A_c$ and $B_c$ takes different forms according to the contact status:
\begin{itemize}
\item First case: sliding contact.

In that case, we have
$$
\mu \|\ta_t\|\geq-\tau_n \qquad \|\ta_t\|\geq \mu \tau_n
$$
then
$$
Proj(\ta,K_\mu)=\tau - \left(\frac{\|\ta_t\|-\mu \tau_n}{1+\mu^2} \right) \left(\frac{\ta_t}{\|\ta_t\|}-\mu \no\right)
$$
and
$$
{\bf Z}_c=\rho(\mu\|\vp^k_t\|\no+\vp^k)+ \left(\frac{\|\ta_t\|-\mu \tau_n}{1+\mu^2} \right) \left(\frac{\ta_t}{\|\ta_t\|}-\mu \no\right)
$$

The computation of the derivatives of ${\bf Z}_c$ provides the matrices $A_c$ and $B_c$:
\begin{itemize}
\item $\displaystyle \frac{\partial {\bf Z}_c}{\partial \sn}  = -\frac{\mu}{1+\mu^2} \left( \frac{\ta_t}{\|\ta_t\|} -\mu \no\right)$
\item $\displaystyle \frac{\partial {\bf Z}_c}{\partial r_{t_1}} = \frac{\tau_{t_1}}{(1+\mu^2)\|\tau_t\|}\left( \frac{\ta_t}{\|\ta_t\|} -\mu \no \right)+
\frac{\|\ta_t\|-\mu\ta_n}{1+\mu^2}\left(\frac{{\bf t}_1}{\|\ta_t\|}-\frac{\tau_{t_1}}{\|\ta_t\|^3}\ta_t \right)$
\item $\displaystyle \frac{\partial {\bf Z}_c}{\partial r_{t_2}} = \frac{\tau_{t_2}}{(1+\mu^2)\|\tau_t\|}\left( \frac{\ta_t}{\|\ta_t\|} -\mu \no \right )+
\frac{\|\ta_t\|-\mu\ta_n}{1+\mu^2}\left(\frac{{\bf t}_2}{\|\ta_t\|}-\frac{\tau_{t_2}}{\|\ta_t\|^3}\ta_t \right)$
\item $\displaystyle \frac{\partial {\bf Z}_c}{\partial v_n} = \rho \no +\frac{\rho\mu}{1+\mu^2}\left(\frac{\ta_t}{\|\ta_t\|}-\mu \no\right)$
\item $\displaystyle \frac{\partial {\bf Z}_c}{\partial v_{t_1}} = \rho\left( {\bf t}_1+\mu\frac{ v_{t_1}}{\|\vp_t\|}\no\right)-\frac{\rho}{1+\mu^2}\left( 
\begin{array}{l}
\left(\frac{\tau_{t_1}}{\|\ta_t\|}-\frac{\mu^2v_{t_1}}{\|\vp_t\|}\right)\left( \frac{\ta_t}{\|\ta_t\|} -\mu \no \right )+\\
\qquad\quad(\|\ta_t\|-\mu \tau_n)\left( \frac{{\bf t}_1}{\|\ta_t\|}-\frac{\tau_{t_1}}{\|\ta_t\|^3}\ta_t\right)
\end{array}\right)$
\item $\displaystyle \frac{\partial {\bf Z}_c}{\partial v_{t_2}} = \rho\left( {\bf t}_2+\mu \frac{v_{t_2}}{\|\vp_t\|}\no\right)-\frac{\rho}{1+\mu^2}\left( 
\begin{array}{l}
\left(\frac{\tau_{t_2}}{\|\ta_t\|}-\frac{\mu^2v_{t_2}}{\|\vp_t\|}\right)\left( \frac{\ta_t}{\|\ta_t\|} -\mu \no \right )+\\
\qquad\quad(\|\ta_t\|-\mu \tau_n)\left( \frac{{\bf t}_2}{\|\ta_t\|}-\frac{\tau_{{\bf t}_2}}{\|\ta_t\|^3}\ta_t\right)\end{array}\right)$
\end{itemize}
\medskip

For a 2D problem, these computations yields:

\begin{itemize}
\item $\displaystyle \frac{\partial {\bf Z}_c}{\partial \sn}  = \frac{\mu}{1+\mu^2} \left( \mu \no- \theta_r {\bf t} \right)$
\item $\displaystyle \frac{\partial {\bf Z}_c}{\partial r_{t}} = \frac{1}{(1+\mu^2)}\left( -\mu \theta_r \no + {\bf t}  \right)$
\item $\displaystyle \frac{\partial {\bf Z}_c}{\partial v_n} = \frac{\rho}{1+\mu^2}\left(\no + \mu \theta_r {\bf t}\right)$
\item $\displaystyle \frac{\partial {\bf Z}_c}{\partial v_{t}} =\frac{\rho\mu}{1+\mu^2}\Bigl((\theta_v+\theta_r)\no + \mu(1-\theta_r\theta_v){\bf t} \Bigr)$
\end{itemize}
\medskip

where $\theta_v=\mbox{sign}(\vp_t)$ and $\theta_r=\mbox{sign}(\ta_t)$.
\bigskip

\item Second case: sticking contact.

In that case, we have
$$
\mu \|\tau_t\|\geq -\tau_n \qquad \|\tau_t\|< \mu \tau_n
$$
then

$$
{\bf Z}_c=\rho (\mu\|\vp^k_t\|\no+\vp^k)
$$
and the computation of the derivatives of ${\bf Z}_c$ reads:
\begin{itemize}
\item $A_c=0_{3\times 3}$
\item $\displaystyle \frac{\partial {\bf Z}_c}{\partial v_n} =\rho \no$
\item$\displaystyle \frac{\partial {\bf Z}_c}{\partial v_{t_1}} = \rho\mu \frac{v_{t_1}}{\|\vp_t\|}\no + \rho {\bf t}_1$
\item$\displaystyle \frac{\partial {\bf Z}_c}{\partial v_{t_2}} = \rho\mu \frac{v_{t_2}}{\|\vp_t\|}\no + \rho {\bf t}_2$
\end{itemize}
\medskip

For a 2D problem, these computations leads to:

\begin{itemize}
\item $A_c=0_{2\times 2}$
\item $\displaystyle \frac{\partial {\bf Z}_c}{\partial v_n} = \rho \no$
\item $\displaystyle \frac{\partial {\bf Z}_c}{\partial v_{t}} =\rho\mu \theta_v \no +\rho {\bf t}$
\end{itemize}
\medskip

\item Third case: no contact.

In that case, the matrices $A_c= \mbox{Id}_{3\times 3}$ and $B_c$ vanishes, and $\chi^{\ell+1}_c=\left\{\begin{array}{c}
0\\ \vp^k\end{array}
\right\}$
\end{itemize}


\subsection{Resolution of the linear system}

Generally, the drawback of a Newton is the computational cost of the linear system to be solved at each iteration.
Here, the particular form of the tangent matrix allows the use of a condensation technique. More precisely,
the linear system to be solved can be written:
\begin{equation}
\left[\begin{array}{cc}
-W & Id_{3\times3} \\
A_c & B_c\\
\end{array}\right]
\left(\begin{array}{c}
\delta \st\\ \delta \vp \end{array}\right) =
\left(\begin{array}{c}
\bf 
-f\\ \bf -g \end{array}\right).
\end{equation}

The  first equation yields $\delta \vp=-{\bf f}+W\delta \st$, and introducing this equality is the second equation leads
to solve the linear system
\begin{equation}
(A_c+B_cW)\delta \st = -{\bf g}+B_c {\bf f}.
\end{equation}
This properties halves the size of the linear system to be solved.

\begin{rmk}
A drawback of the bi-potential framework is that, due to is specificity, it is rather difficult to consider fully coupled problems, where
the contact law and another phenomena, such as electricity or thermic effects are strongly coupled. The other
method presented in this paper has a better property from this point of view because it is based on a  more standard
 mathematical background in the theory of optimization.
\end{rmk}
\subsection{Newton method and enhanced augmented lagrangian method, (SAL) and (EAL)}

In \cite{AC91}, Alart et al. propose another method to solve the frictional contact problem. 
This method has been also used with various improvement
(parallelization, conjugate gradient method for example) to solve multi-contact problems \cite{RA04,RA04b,RDA04,RAD05,RBDA05}.
Even if the coupled frictional contact problem is not an optimization problem anymore, it is always possible to formally formulate a ``quasi"-
optimization problem, for which the constraint set depends on the normal components of the solution as a parameter. The solution is then searched
as a saddle point of a "quasi" augmented Lagrangian of the problem. 

More precisely, 
the global problem on all unknowns that has to be solved at each time step (in place of equation (\ref{locprob1})) has the following form:
\begin{equation}\label{globprob1}
\left\{\begin{array}{l}
{\bf u}={\bf u}^{free} + \WW \st \\
\ \\
 \st\geq 0,\ {\bf u}\geq 0,\  \st\cdot {\bf u}=0.
\end{array}
\right.
\end{equation}

In order to solve this problem, for a given $\st\in \RR^{3\times N_c}$, one can define the cartesian product of infinite half cylinder with section
equal to the ball ${\cal B}(0,\mu r_c)$ of radius $\mu r_c$ by:
\[
{\cal C}(\mu \st)=\prod_{c=1}^{N_c}\RR^+\times {\cal B}(0,\mu r_c)
\]
and then, the granular type frictional contact problem is given by
\begin{equation}\label{AL}
\st\in \mbox{argmin}_{\st\in {\cal C}(\mu \st)}    \frac12 \st \cdot \WW\st +{\bf u}^{free}\cdot \st=  \mbox{argmin}_{\st\in {\cal C}(\mu \st)}  J(\st),
\end{equation}
and the projected gradient method the minimize this problem reads (for each iteration $i$ of the NLGS algorithm):
\begin{equation}\label{gradientAL}
\st^{i+1}=proj(\st^i-\rho({\bf u}^{free}+{\mathbb W}\st^i),\ {\cal C}(\mu{\st^{i+1}})),
\end{equation}
or $\st^{i+1}=proj(\ta^{i+1}, \ {\cal C}(\mu{\st^{i+1}}))$, with $\ta^{i+1}=\st^i-\rho{\bf u}^i$, ${\bf u}^i={\bf u}^{free}+{\mathbb W}\st^i$. This algorithm will be referred to hereinafter as the SAL (Simple Augmented Lagrangian) method.
\medskip

Notice that this method is very closed to the SBP method. More precisely, for a contact $c$, only the descent direction $\up^{c,i}+\mu \|\up^{c,i}_t \|\no$
in (\ref{uzawa1}) is replaced by $\up^{c,i}$ and the projection $\st^{c,i+1}=proj(\ta^{c,i+1},K_{\mu})$ in (\ref{uzawa3}) is replaced by
\[
\left\{\begin{array}{l}
\sn^{c,i+1}=max(0,\tau_n^{c,i+1})\\
\ \\
\st_t^{c,i+1}=\frac{\ta_t^{c,i+1}}{\|\ta_t^{c,i+1}\|}\mu \sn^{c,i+1}.
\end{array}
\right.
\] 
\bigskip

\begin{rmk}
On the contrary, it is possible to see the algorithm developed from the bi-potential formalism as a slight modification
of the algorithm above. Indeed, it is only necessary to change the set ${\cal C}(\st)$ by $\displaystyle {\cal K}=\prod_{c=1}^{N_c}K_\mu$,
and to change the descent direction $\up^{c,i}$ by  $\up^{c,i}+\mu \|\up^{c,i}_t \|\no$ which remains a descent direction for the SAL method, since
$$
\nabla J(\st^{c,i+1})\cdot \Db=-\|\up^{c,i}\|^2-\up^{c,i}\cdot(\mu \|\up^{c,i}_t\|)\no=-\|\up^{c,i}\|^2-\mu u^{c,i}_n\|\up^{c,i}_t\|
$$
which is negative since 
$\mu\in[0,1]$.

\end{rmk}


Then, acting by analogy, we can develop a Newton method to find the minimum of $J$ by seeking the solution as a zero of the function $\tilde{f}(\chi)$ where, for a contact $c$
\[
\tilde{f}_c(\chi)=
\left(\begin{array}{c}
\displaystyle\up^c_{k+1}-\up^{c,free}_{k}-\sum_{\alpha=1}^{N_c}  W_{c\alpha} \st^\alpha \\
\ \\
{\bf \tilde{Z}}^c 
\end{array}\right),
\]
the vector ${\bf \tilde{Z}}^c$ is the error on the prediction of the reaction
\begin{equation}\label{predbipo}
{\bf \tilde{Z}}^c(\st^c,\up^c_{k+1})=\st^c-proj(\ta^c_{k+1},{\cal C}_c(\mu\ta^c_{k+1})),
\end{equation} 

and the set ${\cal C}_c(\mu\st^c)$ is the set of admissible forces ${\cal C}_c(\mu\st^c)=\RR^+\times {\cal B}(0,r_c)$. 
This method will be refered as the EAL (Enhanced Augmented Lagrangian) method hereafter.
\medskip

Then, as bellow, we have
three cases in the computation of the tangent matrix $\left[\frac{\partial \tilde{f}}{\partial \chi^c}(\chi^\ell)\right]$:

\begin{itemize}
\item First case: sliding contact ($\tau_n>0,\ \ta_t\geq \mu \tau_n$)

We have: $proj(\ta^c,{\cal C}_c(\mu\ta^c))=\tau_n\no + \frac{\ta_t}{\|\ta_t\|}\mu\tau_n {\bf t}$ and
$ \tilde{\bf Z}_c=\rho v_n\no-\frac{\ta_t}{\|\ta_t\|}\mu\tau_n +\st_t$.

The computation of the derivatives of $ \tilde{\bf Z}_c$ provides the matrices $A_c$ and $B_c$:
\begin{itemize}
\item $\displaystyle \frac{\partial \tilde{\bf Z}_c}{\partial \sn}  = -\mu \frac{\ta_t}{\|\ta_t\|}$
\item $\displaystyle \frac{\partial  \tilde{\bf Z}_c}{\partial r_{t_1}} = {\bf t}_1-\mu \tau_n\left( \frac{{\bf t}_1}{\|\ta_t\|}-\frac{\ta_{t_1}}{\|\ta_t\|^3}\ta_t \right)$
\item $\displaystyle \frac{\partial  \tilde{\bf Z}_c}{\partial r_{t_2}} = {\bf t}_2-\mu \tau_n\left( \frac{{\bf t}_2}{\|\ta_t\|}-\frac{\ta_{t_2}}{\|\ta_t\|^3}\ta_t \right) $
\item $\displaystyle \frac{\partial  \tilde{\bf Z}_c}{\partial v_n} = \rho \left(\no +\mu \frac{\ta_t}{\|\ta_t\|}\right)$
\item $\displaystyle \frac{\partial  \tilde{\bf Z}_c}{\partial v_{t_1}} = -\rho \mu \tau_n \left( \frac{{\bf t}_1}{\|\ta_t\|}-\frac{\ta_{t_1}}{\|\ta_t\|^3}\ta_t \right) $
\item $\displaystyle \frac{\partial  \tilde{\bf Z}_c}{\partial v_{t_2}} = -\rho \mu \tau_n \left( \frac{{\bf t}_2}{\|\ta_t\|}-\frac{\ta_{t_2}}{\|\ta_t\|^3}\ta_t \right)$
\end{itemize}
\medskip

For a two dimensional problem, these computations yields
$$
A_c=\left(\begin{array}{cc}
0 & 0 \\
-\mu \theta_r & 1\\
\end{array}\right)\qquad
B_c=\left(\begin{array}{cc}
1 & 0\\
\mu \theta_r & 0\\
\end{array}\right).
$$

\item Second case: sticking contact ($\tau_n>0$, $\ta_t<\mu \tau_n$)

$proj(\ta^c,{\cal C}_c(\mu\ta^c))=\ta^c$
and the computation of the derivatives of ${\bf Z}_c$ reads
\begin{itemize}
\item $A_c=0_{3\times 3}$
\item $B_c=\rho \mbox{Id}_{3\times 3}$
\end{itemize}
\medskip

\item Third case: no contact ($\ta_n\leq0$)

$proj(\ta^c,{\cal C}_c(\mu\ta^c))=0$, then the matrices $A_c= \mbox{Id}_{3\times 3}$ and $B_c$ vanishes, and $\chi^{\ell+1}_c=\left\{\begin{array}{c}
0\\ \vp^k\end{array}
\right\}$

\end{itemize}

\subsection{The global stopping (convergence) criterion}


We present in this paragraph the convergence criterion on the global non linear Gauss-Seidel iterations. This
criterion, developed from that proposed in \cite{FHdS02} has been extended in the case of the Newton
and bi-potential (EBP) method, where some term are naturally vanishing in the original Uzawa and bi-potential (SBP) method.
This criterion $\varepsilon_{glob}$ has been written in such a way that if the solution verify that $\varepsilon_{glob}$
is sufficiently small, then this solution has good properties on the equation of motion and Signorini Coulomb contact law.
Consequently, this criterion stays valid for the methods developed with the augmented lagrangian (SAL and EAL methods).
\medskip

This criterion can be stated:
\begin{equation}\label{estimator}
\varepsilon_{glob}=\frac{1}{N_c}\sum_{c=1}^{N_c}\left[\varepsilon^c_{motion}+ \varepsilon^c_{proj}+\varepsilon_{b_c}+\varepsilon^c_{pen}
\right]
\end{equation}
\goodbreak

where:
\begin{itemize}
\item $\varepsilon^c_{motion}=\|\up^c-\up_m^c\|$ where $\up^c_m=\up^{c,i}+\sum_{\alpha=1}^{N_c} W_{c\alpha} \st^\alpha$, so $\varepsilon_{motion}$ measures the error on the equation of motion (see equation  (\ref{locprob1}), this term vanishes for the SBP and SAL method);
\item $\varepsilon^c_{proj}=\sqrt{\|\st^c-proj(\st^c,{K_\mu})\|^2}$ is the error for the projection on the Coulomb cone (vanishing for the SBP method);
\item $\varepsilon_{b_c}=\Bigl|\up^c\cdot \st^c+\mu \sn^c \|\up^c_t\| \Bigr|$ is the absolute value of the bi-potential that has to vanish if and only if the couple
$(\up^c,\st^c)$ verifies the Signorini Coulomb contact law (see formula \ref{SCrules});
\item $\varepsilon^c_{pen}=-\min(0,\un^c)$ is the value of the penetration.

\end{itemize}
\bigskip

\begin{rmk}
One can notice that is absolutely necessary to verify in the criterion that there is no penetration, because
nothing in the presented algorithm ensures that is condition is verify at the end of the loop.
Moreover, if this condition is not satisfied, the rest  of bi-potential can be negative or equal to zero, even if the couple
$(\up,\st)$ is not a solution.
\end{rmk}

\section{Numerical results}

We present in the section three numerical examples with an increasing complexity.

In these computations, the descent parameter $\rho$ is taken in such a way that the result is optimal, in terms of time computing.
Denoting $\bar \rho= \frac{m_im_j}{m_i+m_j}\frac1{\Delta t}$, for the SBP and the SAL methods, we have chosen
$\rho=0.6 \bar \rho$, 
whereas for the EBP and the EAL methods, it is better the take $\rho=\bar \rho$. 
We recall that it has been show that, for the the bi-potential method 
(see for example \cite{FJCM05}) and
the augmented lagrangian method (see for example \cite{RA04}), the parameter $\rho$ has to verify $\rho< 2\bar \rho$ in order to ensure the convergence.
Generally, for these two methods, the convergence is very sensitive on this parameter.
We will show in the last paragraph of this study that for the EBP method, the parameter $\rho$ can be taken
in a large range around the value $\bar \rho$ without changing dramatically the convergence
of the method.

At each iteration of the NLGS algorithm,  the Newton algorithm is stopped either if the convergence is obtained ($\varepsilon_{Newt}^c\leq 10^{-5}$), or if
the number of iteration of the Newton algorithm reached 100 when there is no convergence.

\subsection{Ball sliding on a plane}

In this first example, we consider a ball placed on a table with an initial horizontal velocity equal to 1.5 m$\cdot$ s$^{-1}$. The ray of the ball is equal to $5\cdot 10^{-3}$ m,
and the friction coefficient between wall and ball is equal to $\mu=0.7$. The time step of discretization is equal to $10^{-4}$ s. In this experiment, the ball first slides
on the table, and then the ball rolls without sliding. The global stoping criterion is equal to $\varepsilon_{glob}=10^{-10}$.

\begin{figure}[htbp]
\begin{center}
\includegraphics[width=6.0cm]{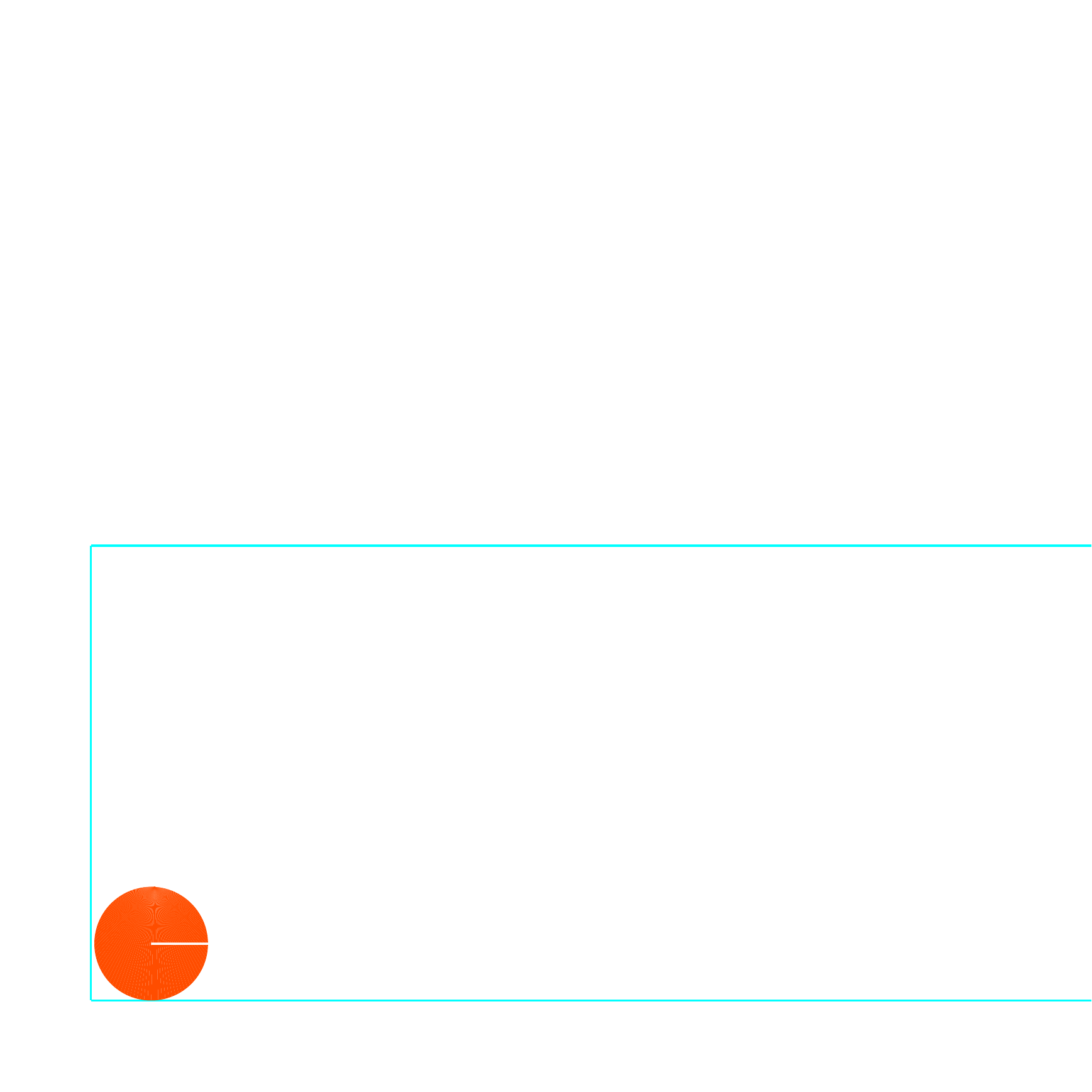}\hspace{0.3cm}\includegraphics[width=6.0cm]{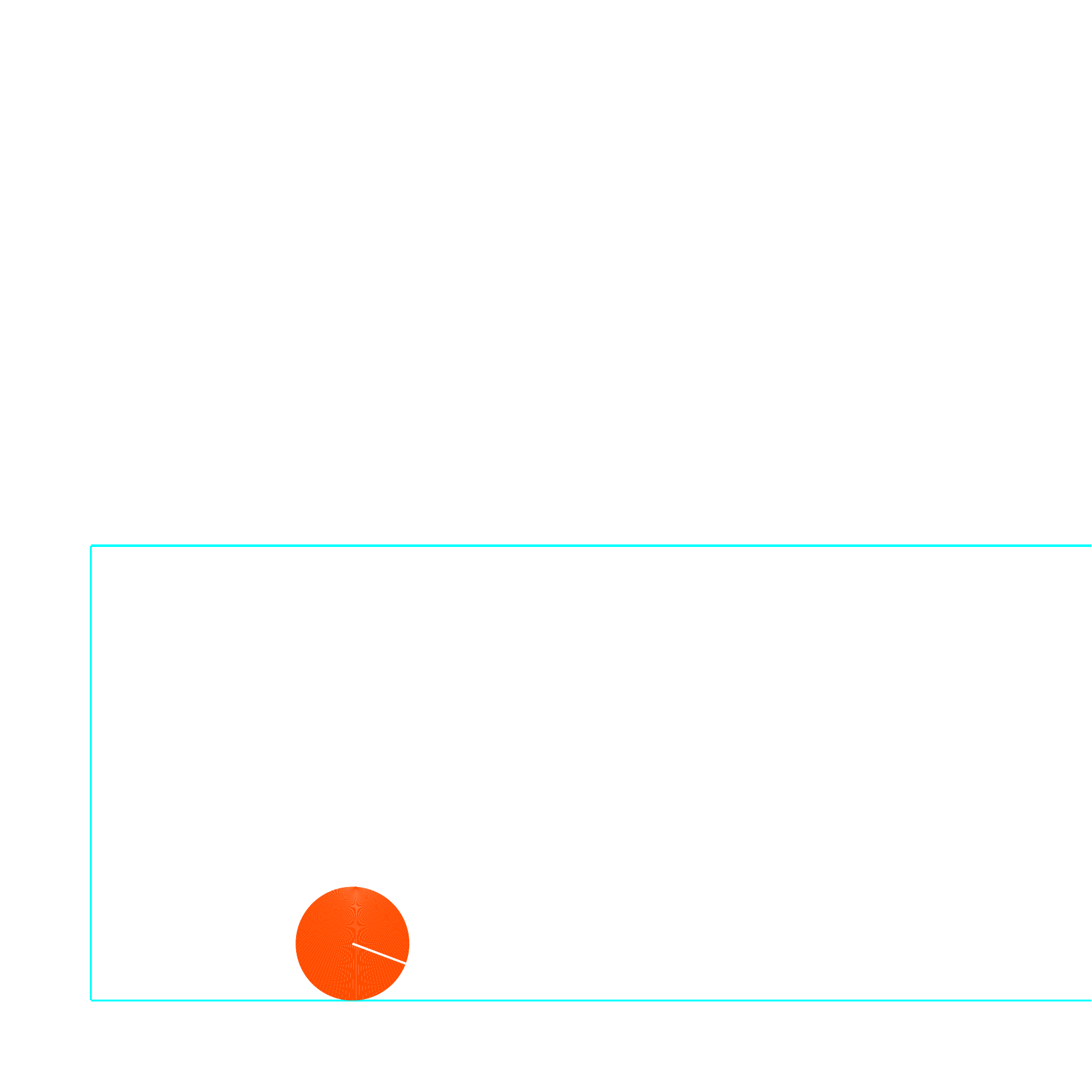}
\caption{Example 1 -- A ball is launched with an initial horizontal velocity (left). First, the ball slides. Then, the ball rolls without slipping (right).}
\label{fig1}
\end{center}
\end{figure}

\begin{table}[htbp]
\begin{center}
\begin{tabular}{|l|c|c|c|}
\hline
Method & Number of &Error & Total \\
 & NLGS iterations & $\varepsilon_{glob}$& CPU time (s)\\
  & (last time step)&  (last time step) & \\

\hline
SBP&18&$0.384\cdot10^{-10}$&9.44\\
SAL&18&$0.384\cdot10^{-10}$&9.28\\
EBP&1&0&8.83\\
EAL & 1&$0.175\cdot10^{-13}$ & 8.78 \\
\hline
\end{tabular}
\caption{Comparison of the results obtained by the four methods on the first example (after the 2000th time step).}
\label{Tab:tab1}
\end{center}
\end{table}%

\begin{figure}[htbp]
\begin{center}
\includegraphics[width=12.0cm]{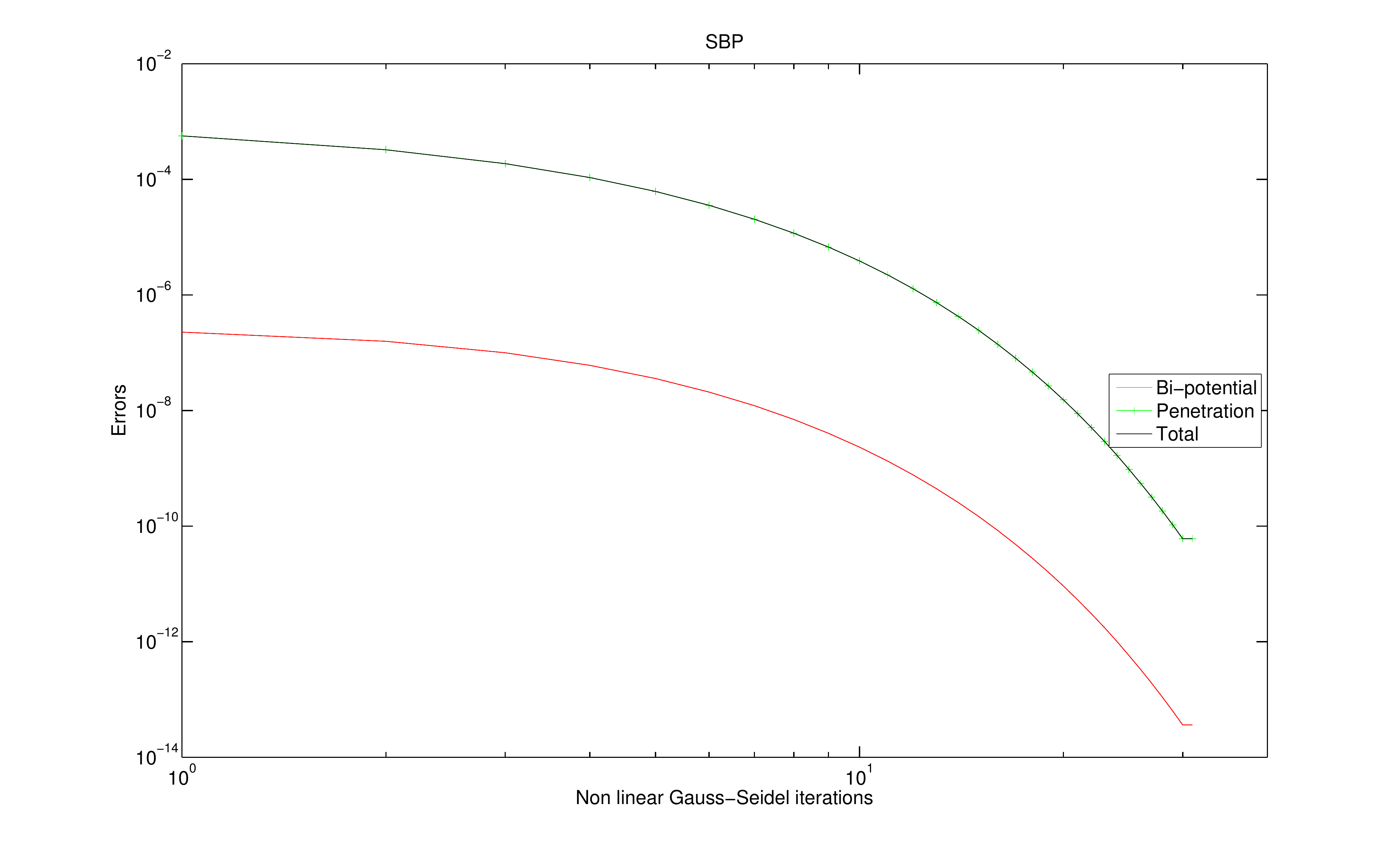}
\caption{Example 1 -- Convergence for the standard bi-potential based method, 5$^{th}$ iteration}
\label{conv1-1}
\end{center}
\end{figure}

\begin{figure}[htbp]
\begin{center}
\includegraphics[width=12.0cm]{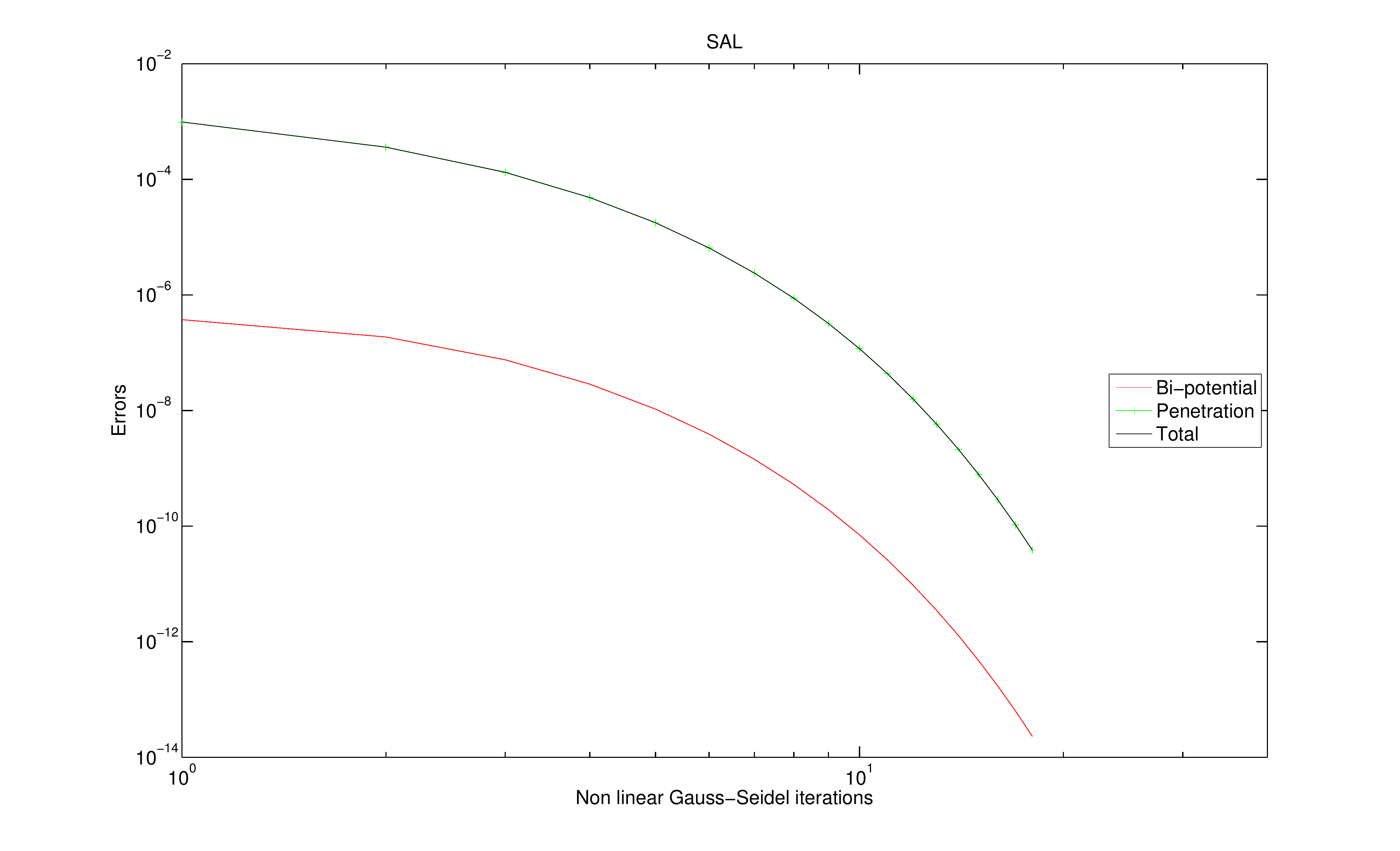}
\caption{Example 1 -- Convergence for the standard augmented lagrangian method, 5$^{th}$ iteration}
\label{conv1-2}
\end{center}
\end{figure}

\begin{figure}[htbp]
\begin{center}
\includegraphics[width=12.0cm]{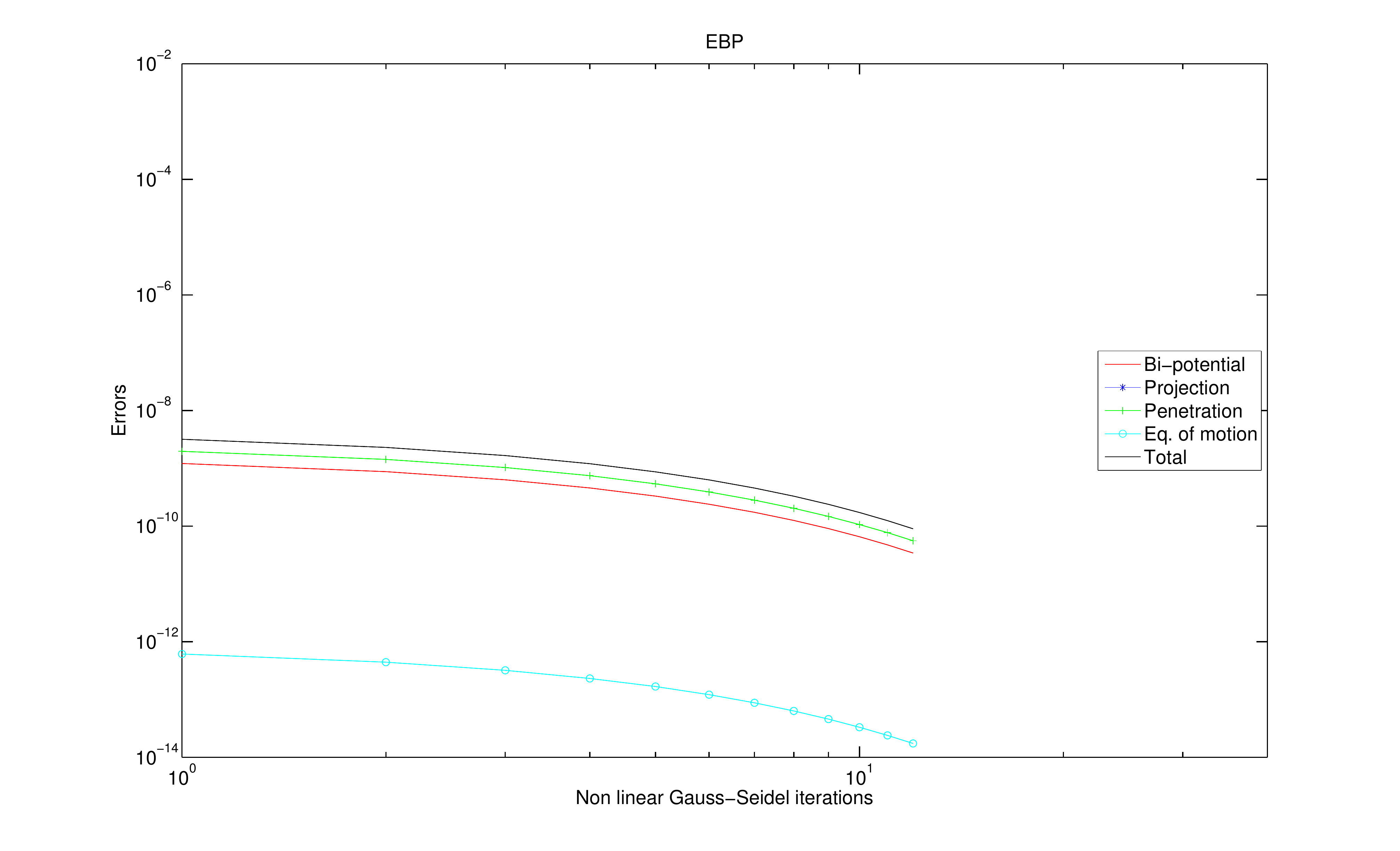}
\caption{Example 1 -- Convergence for the Newton and bi-potential method, 5$^{th}$ iteration}
\label{conv1-3}
\end{center}
\end{figure}

We can observe from these numerical results that the error coming from the projection is very small for the four
methods. The Standard Bi-Potential (SBP) method and the Standard Augmented Lagrangian method (SAL) give
very closed results, both in term of quality (see figures \ref{conv1-1} and \ref{conv1-2}) and in term 
of time computing (see table \ref{Tab:tab1}). Nevertheless, we can notice that the time computing
is smaller with the SAL method, because there is less computations at each iteration (no term such as $\|\up_t\|$ and
projection easier to compute for example). The Enhanced Bi-Potential method provides better results, both in term
of quality (see figure \ref{conv1-3}) and in term of time computing (6.5\% better). The Enhanced Augmented Lagrangian method converges 
after the first Non Linear Gauss Seidel iteration for every time steps, and consequently, this is the faster method on this example
(7\% faster than the SBP method).

\subsection{Sedimentation of 4 balls in a box}

In this second experiment, we consider the
sedimentation of 4 balls of radius ranging from $4\cdot10^{-4}$ m to  $5\cdot10^{-4}$ m.
For the computations, the time step of discretization is equal to 
$\Delta t= 10^{-4}$ s., and the Non linear Gauss-Seidel loop is stopped either if the
the global stopping criterion on the NLGS method is equal to $\varepsilon_{glob}=10^{-10}$,
or after 5000 iterations if there is no convergence (this case never occurs in this experiment). 
The friction coefficient between the balls and between the
balls and the walls is equal to $\mu=0.3$.

\begin{figure}[htbp]
\begin{center}
\includegraphics[width=6.0cm]{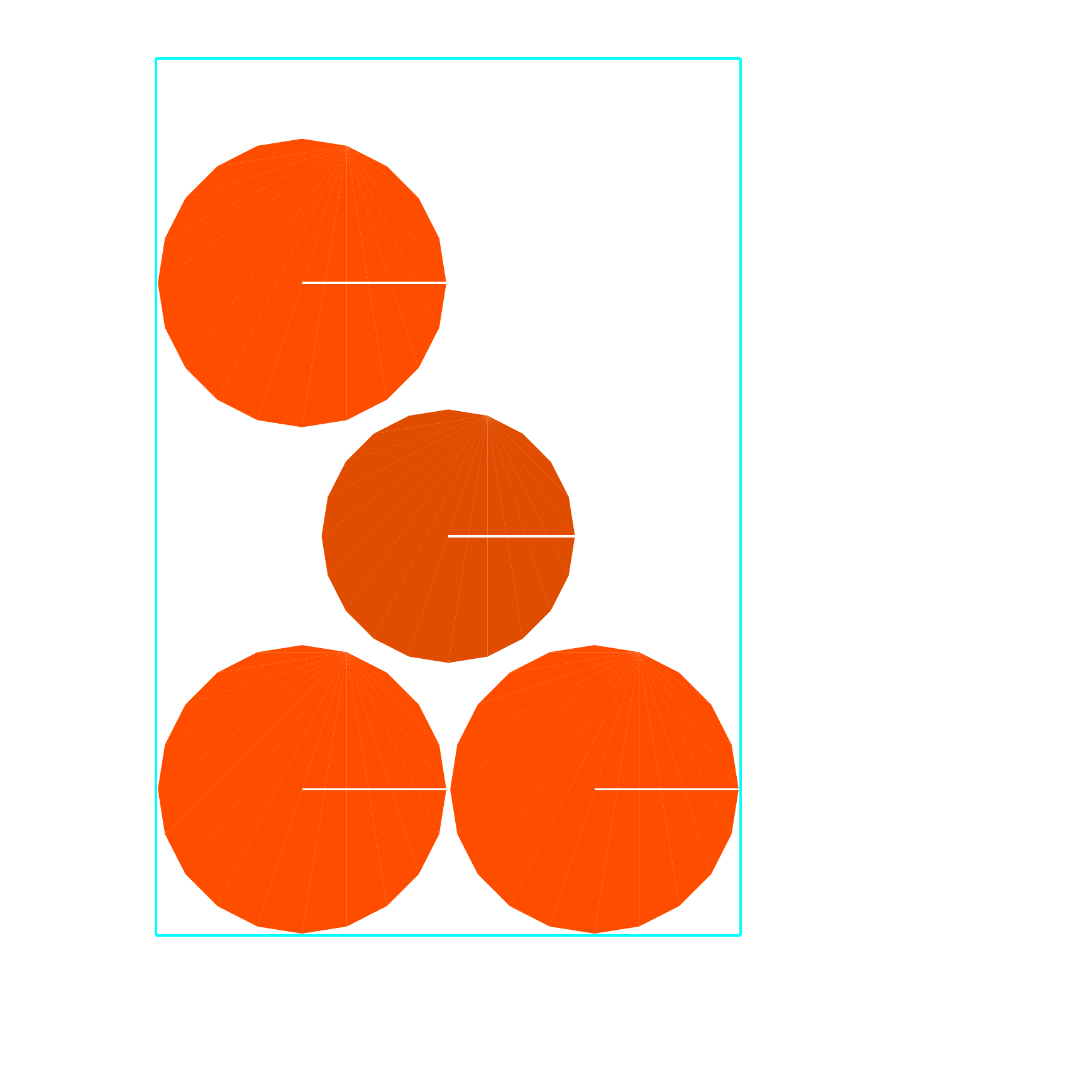}\hspace{0.3cm}\includegraphics[width=6.0cm]{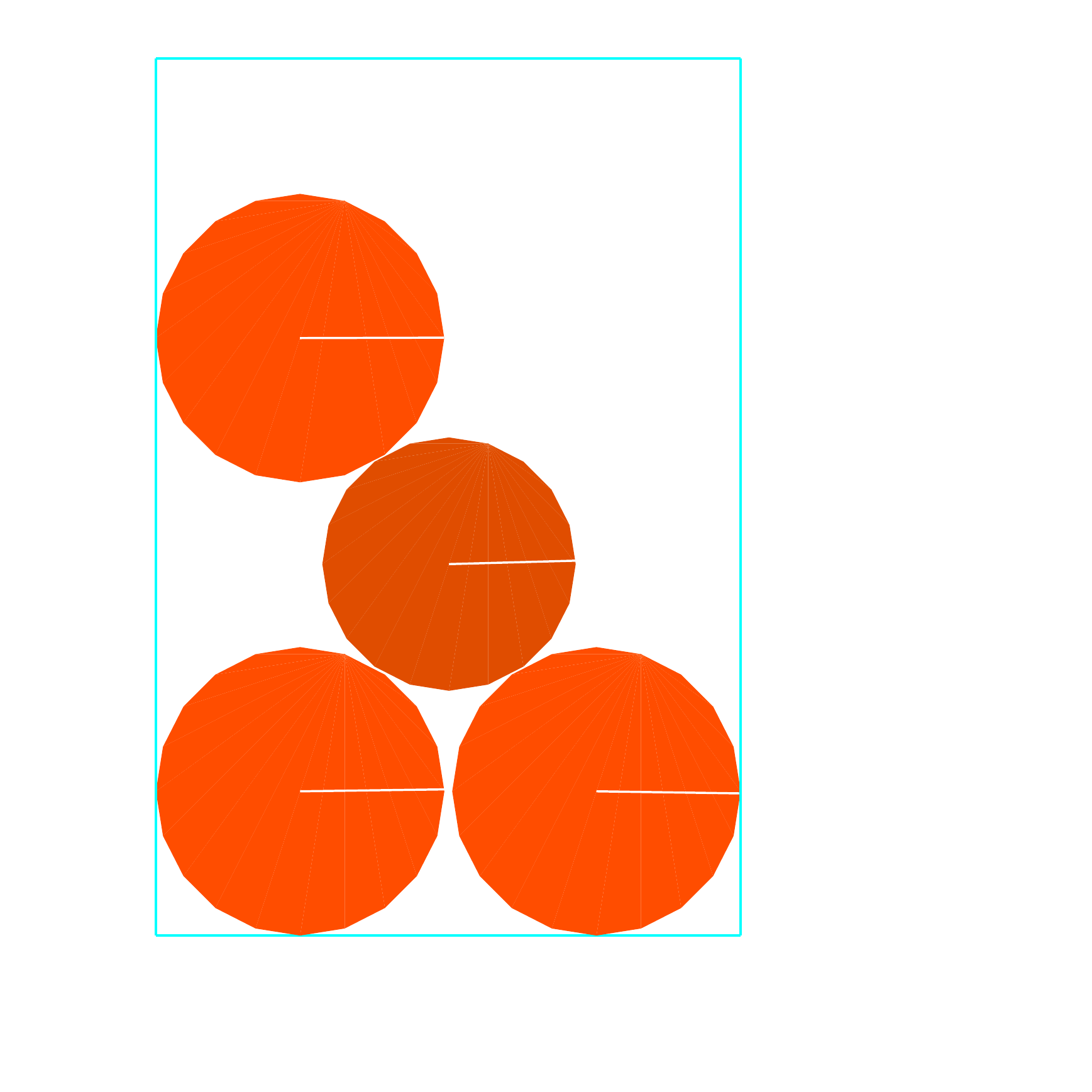}
\caption{Example 2 -- Sedimentation of four balls under the gravity effect.}
\label{fig2}
\end{center}
\end{figure}

\begin{table}[htbp]
\begin{center}
\begin{tabular}{|l|c|c|c|c|}
\hline
Method &Number of  &Error & Maximal &  Total\\
 &NLGS iterations& $\varepsilon_{glob}$ &  penetration & CPU time (s)\\
  & (last time step)&  (last time step)&(last time step)& \\
\hline
SBP&305&$0.949\cdot10^{-12}$ &$0.310\cdot10^{-11}$&2.92\\
SAL&301&$0.980\cdot10^{-12}$&$0.340\cdot10^{-11}$&2.87\\
EBP&161&$0.635\cdot10^{-12}$&$0.641\cdot10^{-12}$&2.59\\
EAL & 158& $0.973\cdot10^{-12}$& $0.208\cdot10^{-19}$& 2.43\\
\hline
\end{tabular}
\caption{Comparison of the results obtained by the four methods on the second example (after the 1000th time step)\label{Tab:tab2}}
\end{center}

\end{table}%

\begin{figure}[htbp]
\begin{center}
\includegraphics[width=12.0cm]{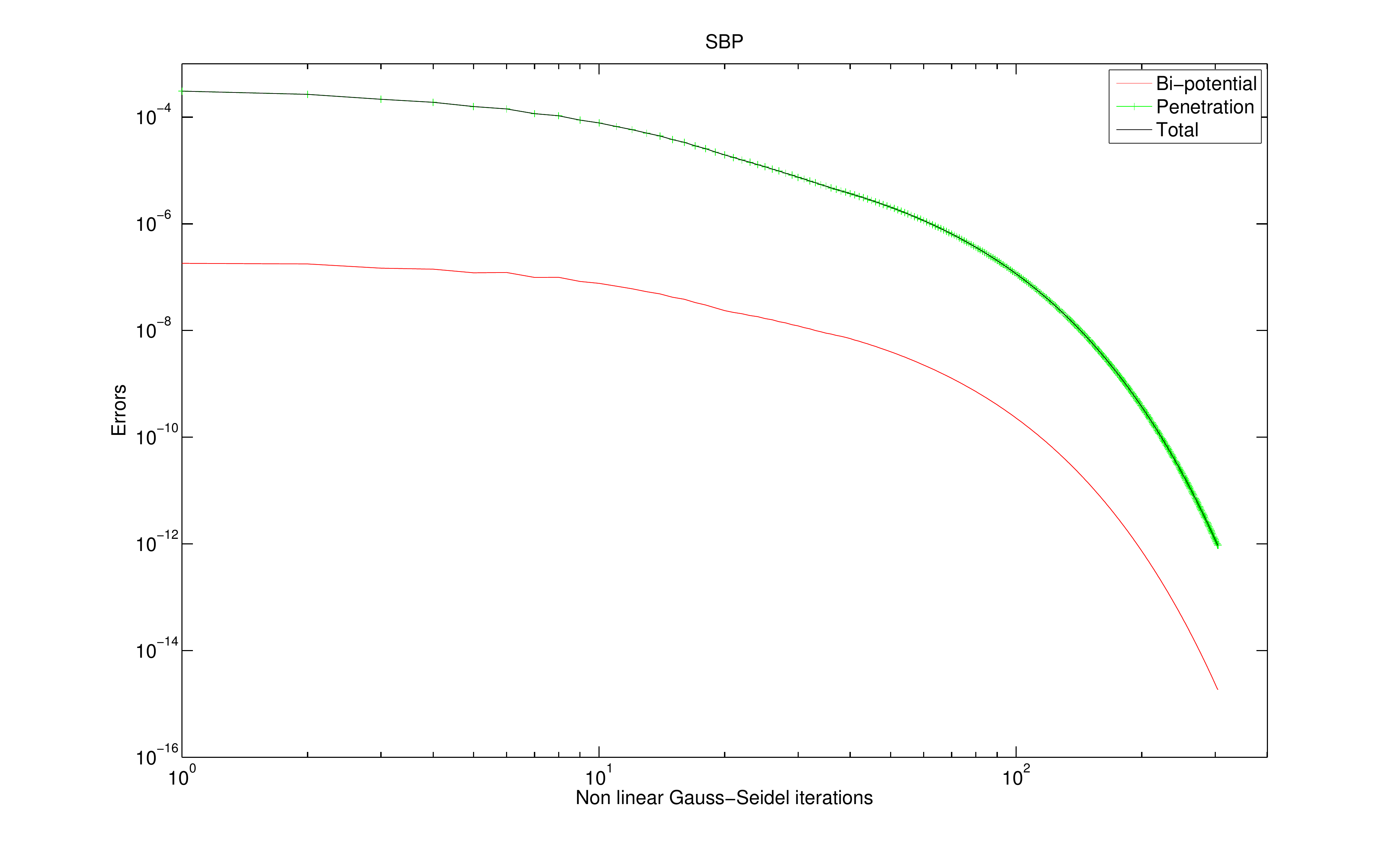}
\caption{Example 2 -- Convergence of the non-linear Gauss-Seidel iterations for the standard bi-potential based method (1000th time step)
The two last curves overlaps, showing that the global error is governed by the error of penetration.}
\label{conv2-1}
\end{center}
\end{figure}

\begin{figure}[htbp]
\begin{center}
\includegraphics[width=12.0cm]{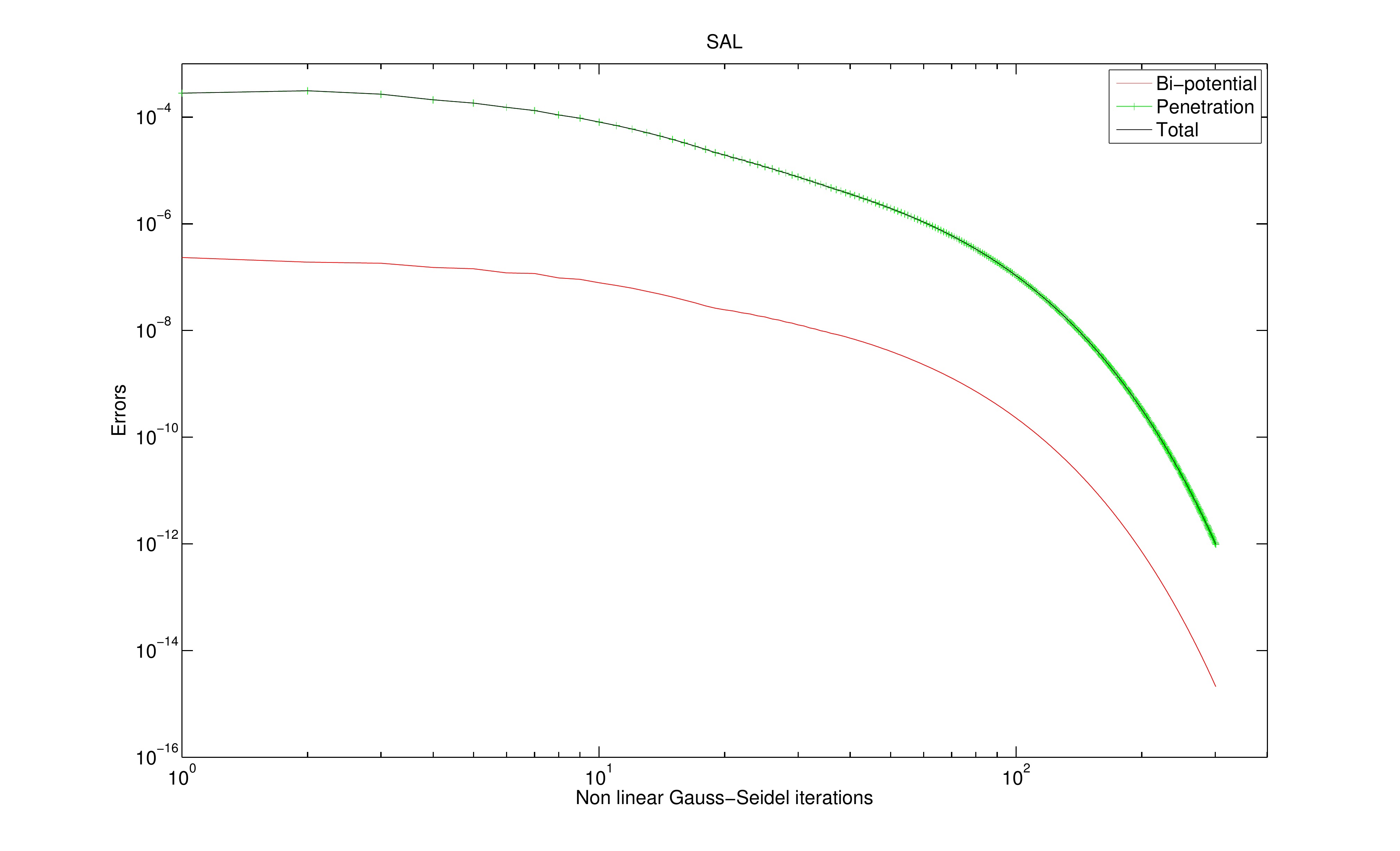}
\caption{Example 2 -- Convergence of the non-linear Gauss-Seidel iterations for the standard augmented lagrangian method (1000th time step).
The two last curves overlaps, showing that the global error is governed by the error of penetration.}
\label{conv2-2}
\end{center}
\end{figure}

\begin{figure}[htbp]
\begin{center}
\includegraphics[width=12.0cm]{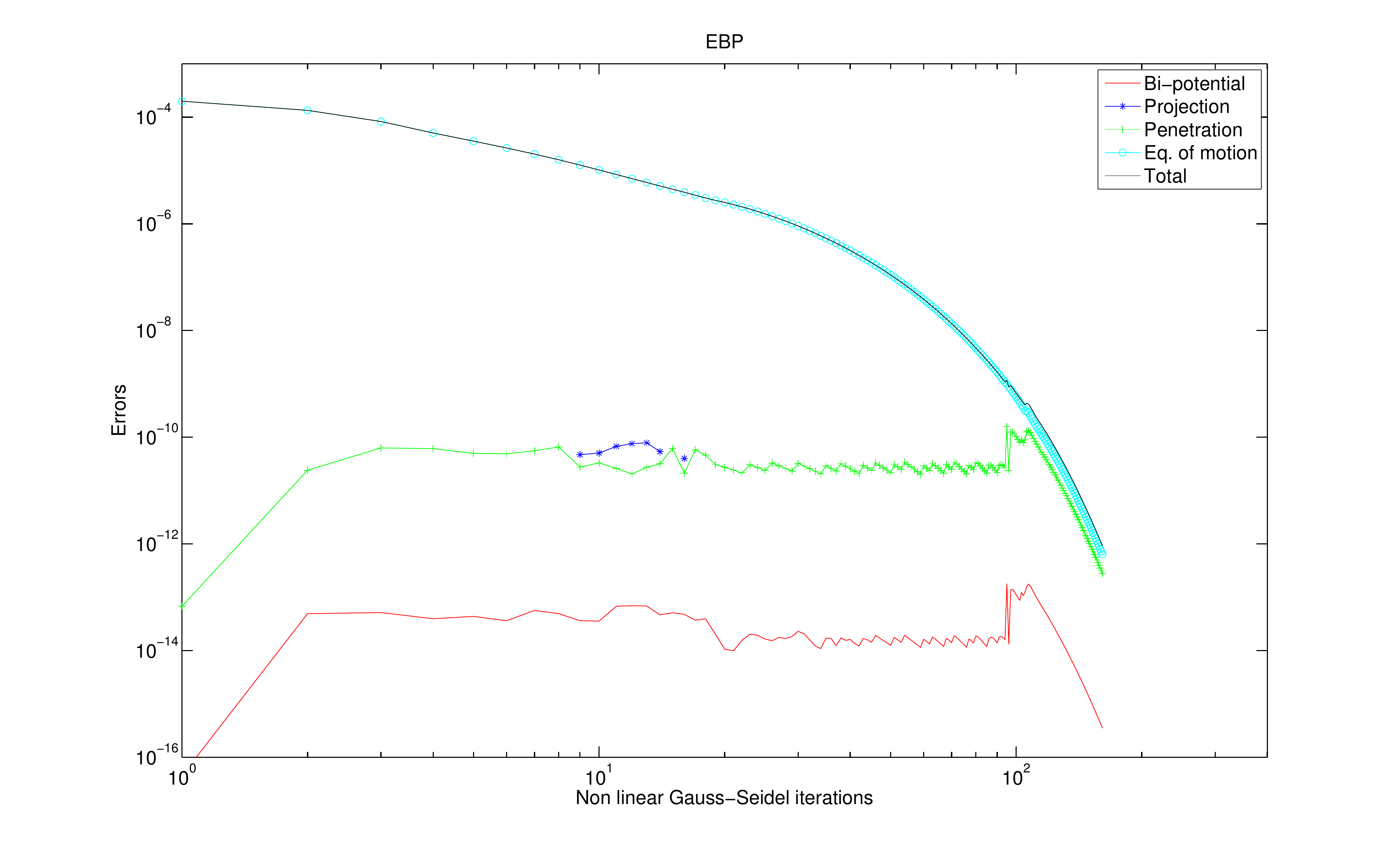}
\caption{Example 2 -- Convergence of the non-linear Gauss-Seidel iterations for the Newton and bi-potential method (1000th time step).
The two last curves overlaps, that shows that the global error is governed by the error on the equations of motion.}
\label{conv2-3}
\end{center}
\end{figure}

\begin{figure}[htbp]
\begin{center}
\includegraphics[width=12.0cm]{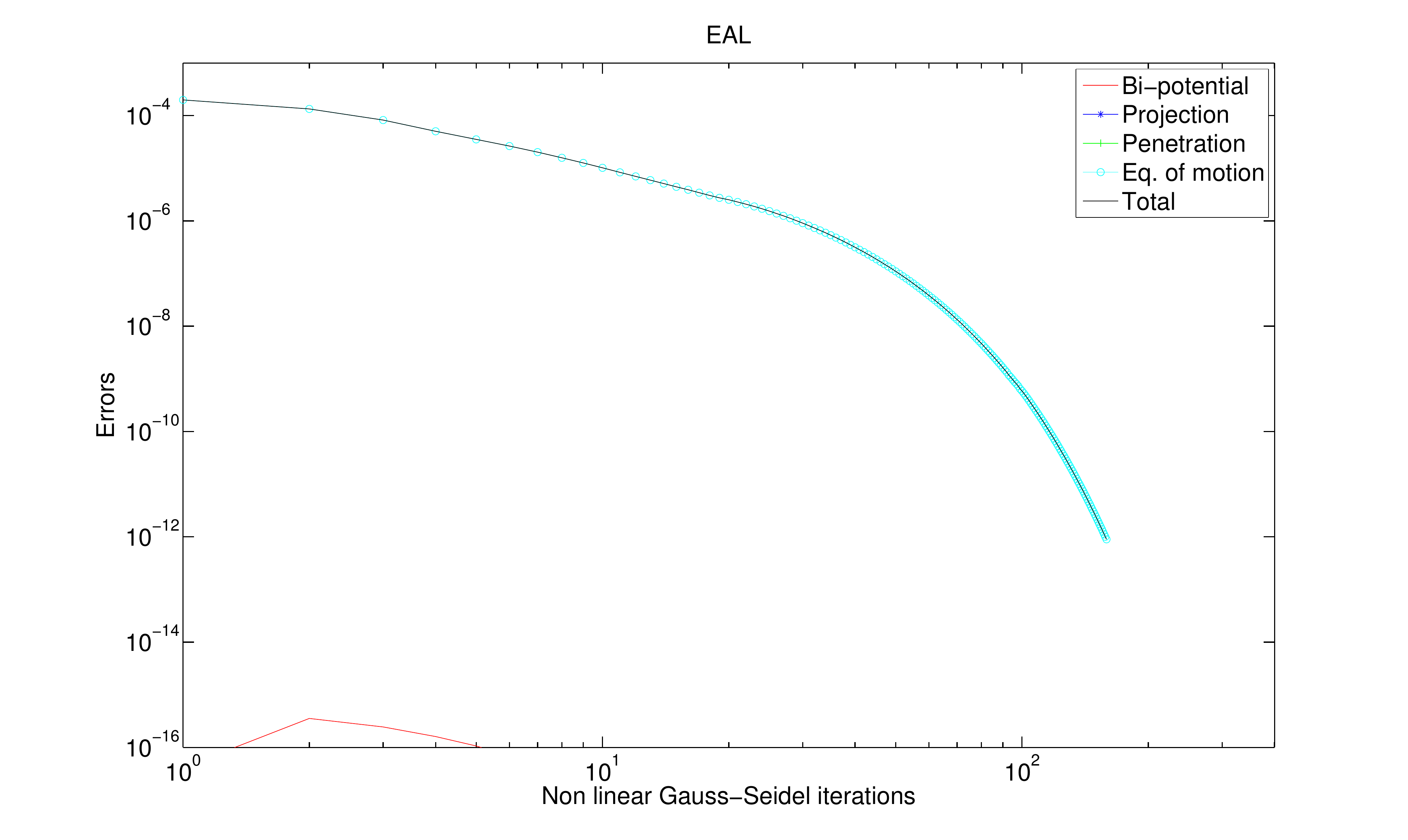}
\caption{Example 2 -- Convergence of the non-linear Gauss-Seidel iterations for the Newton and Augmented Lagrangian method (1000th time step).
The two last curves overlaps, and the other ones does not appear on the figure because the corresponding errors are lower than $10^{-16}$.}
\label{conv2-4}
\end{center}
\end{figure}

Like in the previous simulation,, the SBP and the SAL method methods provide very similar results (see figures \ref{conv2-1} and \ref{conv2-2}). 
For these two methods, we can notice that here 
the global error is essentially due to the penetrations. The SAL method is 2\% faster than the SBP method (table \ref{Tab:tab2}).

Results obtained by the EBP method are better (figure \ref{conv2-3}), and here the overall error is governed by the error on the
equations of the motion. The EBP method is 11\% faster than the SBP method, and the penetration is 5 times smaller.
In this example the EAL is the faster method (16,8\% faster than the SBP method), and the penetration is very small (see figure \ref{conv2-4}).

\subsection{Sedimentation of 500 balls}
\begin{figure}[htbp]
\begin{center}
\includegraphics[width=6.0cm]{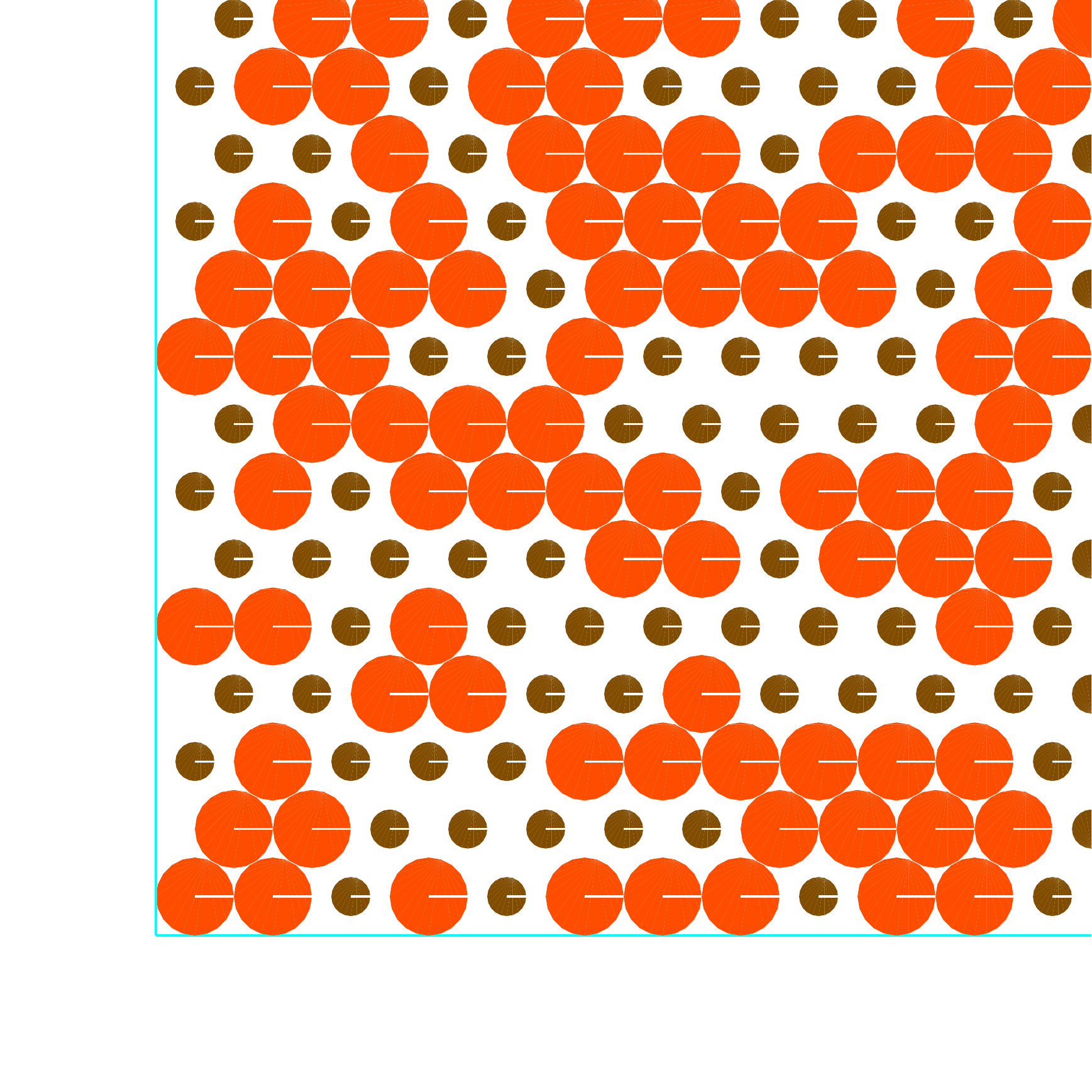}\hspace{0.3cm}\includegraphics[width=6.0cm]{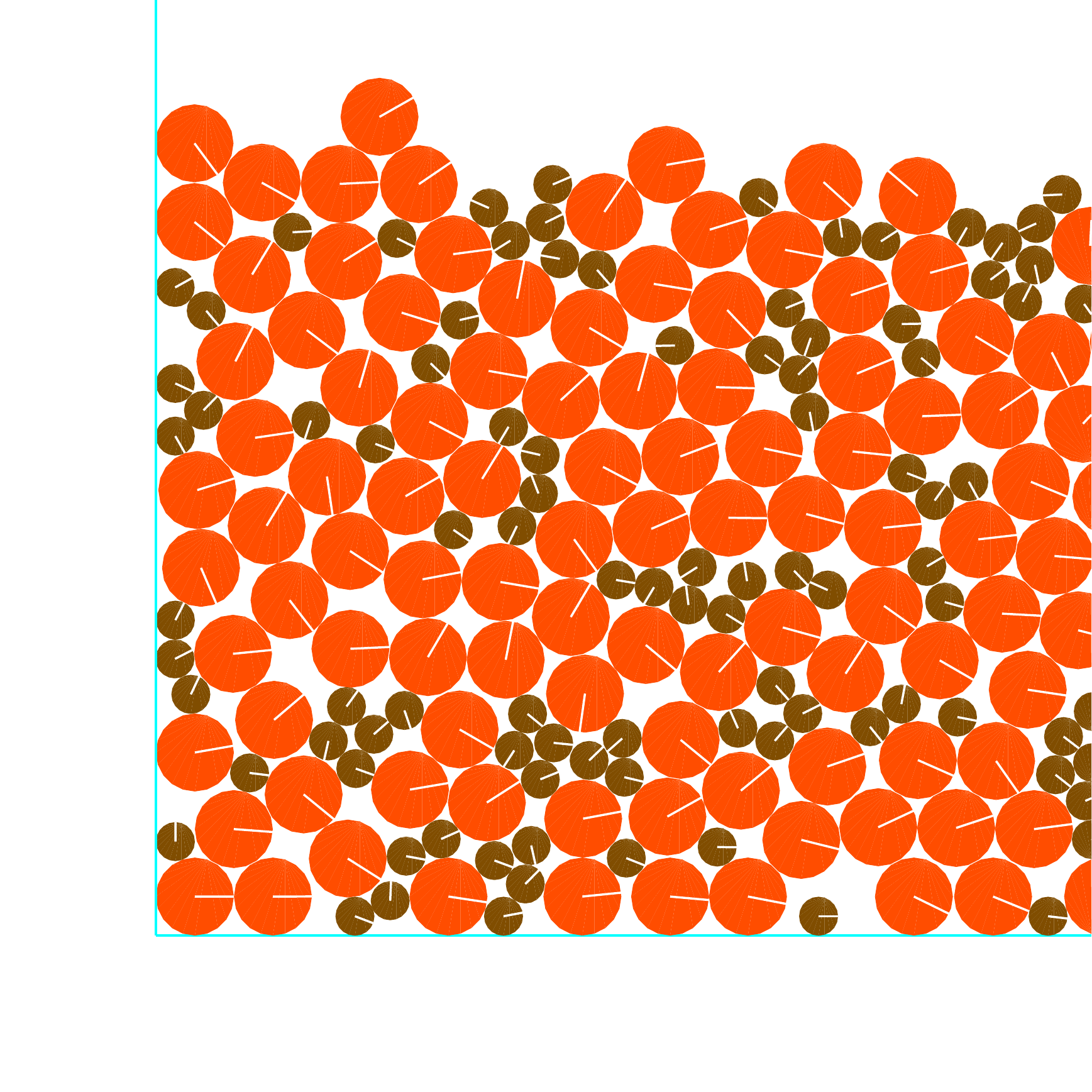}
\caption{Example 3 -- Zoom on balls falling under the gravity effect. Initial configuration on the left, final configuration on the right.}
\label{fig3}
\end{center}
\end{figure}

In this example, we consider the
sedimentation of 500 balls (see figure \ref{fig3}) of radii ranging from $2.5\cdot10^{-4}$ m to  $5\cdot10^{-4}$ m, 
the time step of discretization is equal to $\Delta t=5\cdot 10^{-5}$ s, and the Non linear Gauss-Seidel loop is stopped if the
global estimator (\ref{estimator}) verifies $\varepsilon_{glob}\leq 10^{-12}$ or after after 5000 iterations if there is no convergence.
The friction coefficient between the balls and between the balls and the walls is equal to $\mu=0.3$.

The results in table \ref{Tab:tab3-1} are obtained after 1000 time steps.


\begin{table}[htbp]
\begin{center}
\begin{tabular}{|l|c|c|c|c|}
\hline
Method & Number of &Error& Maximal   & Total\\
& NLGS iterations& $\varepsilon_{glob}$&  penetration  & CPU time (s)\\
 & (last time step)&  (last time step)&(last time step)& \\
\hline
SBP&5000&$0.119\cdot 10^{-6}$&$0.213\cdot 10^{-5}$ & 1092.95\\
SAL&5000&$0.135\cdot 10^{-6}$&$0.533\cdot 10^{-5}$ & 973.31\\
EBP&5000&$0.156\cdot 10^{-6}$&$0.286\cdot 10^{-6}$ & 854.31\\
EAL & 5000& $0.101\cdot 10^{-6}$& $0.390\cdot 10^{-17}$& 916.65\\
\hline
\end{tabular}
\caption{Comparison of the results obtained by the four methods on the third example (after the 1000th iteration, $N_{max}=5000$ iterations)\label{Tab:tab3-1}}
\end{center}
\end{table}

\begin{figure}[htbp]
\begin{center}
\includegraphics[width=12.0cm]{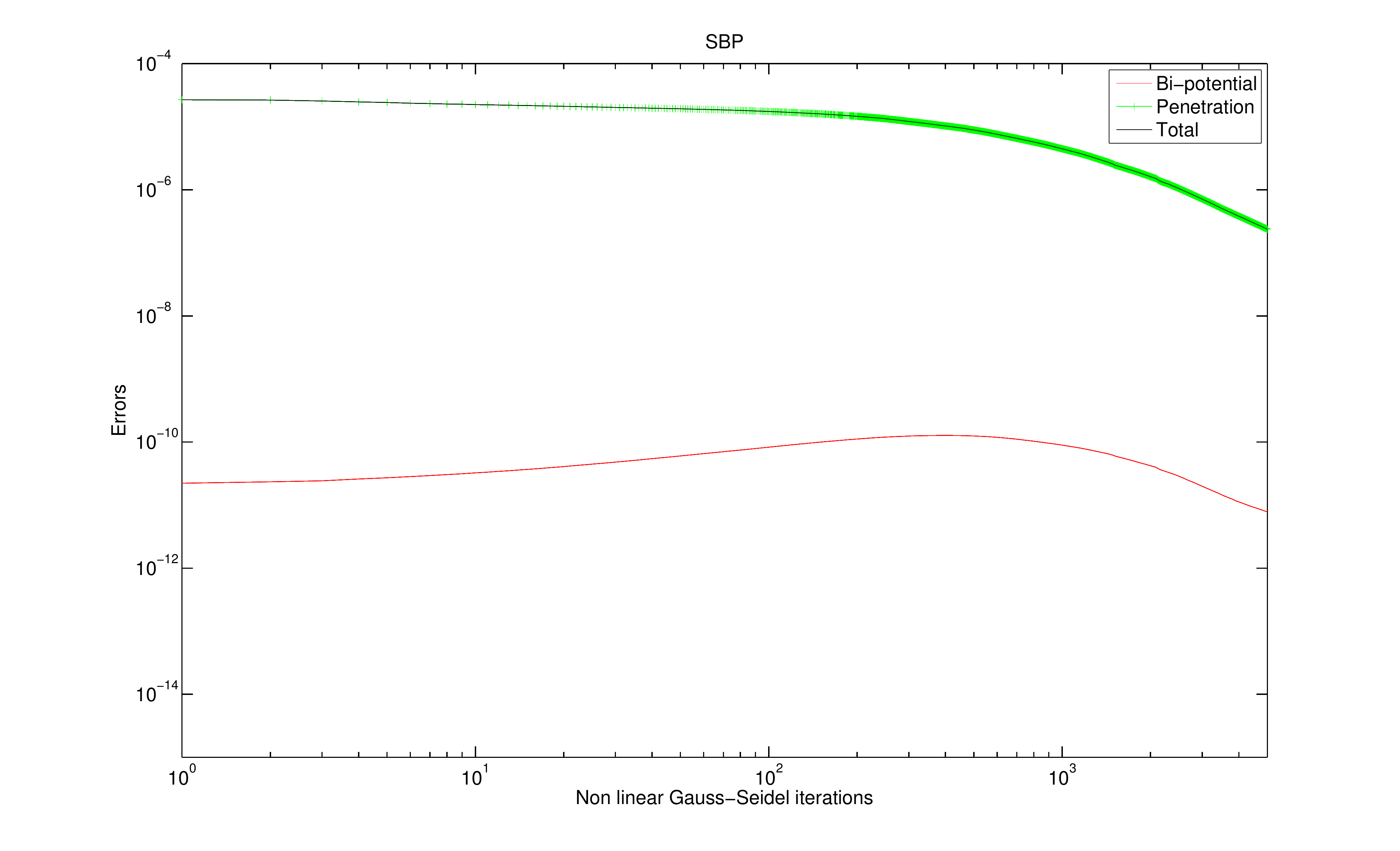}
\caption{Example 3 -- Convergence of the non-linear Gauss-Seidel iterations for the standard bi-potential based method  (1000th time step). The two last curves collapse.}
\label{conv3-1}
\end{center}
\end{figure}

\begin{figure}[htbp]
\begin{center}
\includegraphics[width=12.0cm]{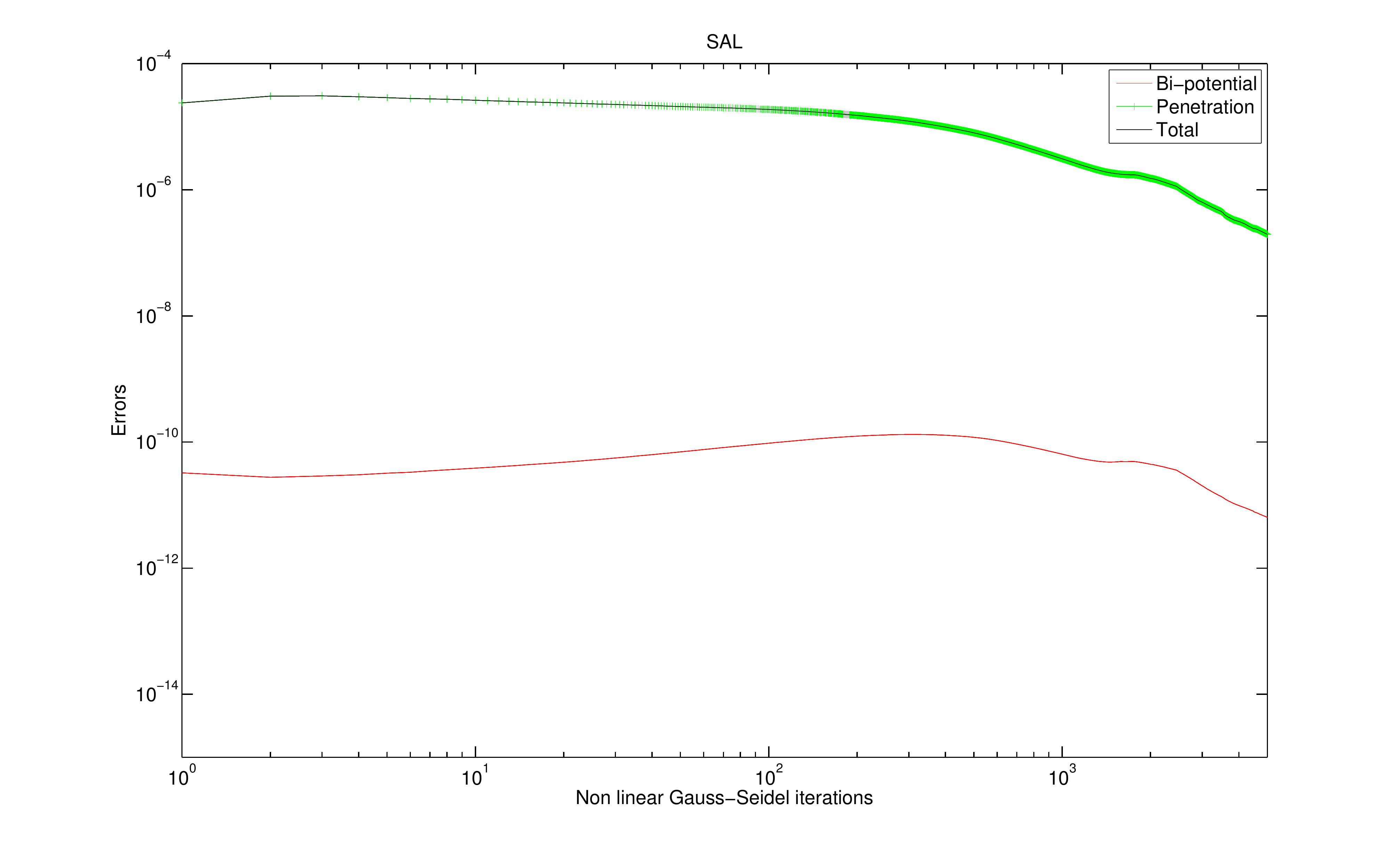}
\caption{Example 3 -- Convergence of the non-linear Gauss-Seidel iterations for the standard augmented lagrangian method  (1000th time step). The two last curves collapse.}
\label{conv3-2}
\end{center}
\end{figure}

\begin{figure}[htbp]
\begin{center}
\includegraphics[width=12.0cm]{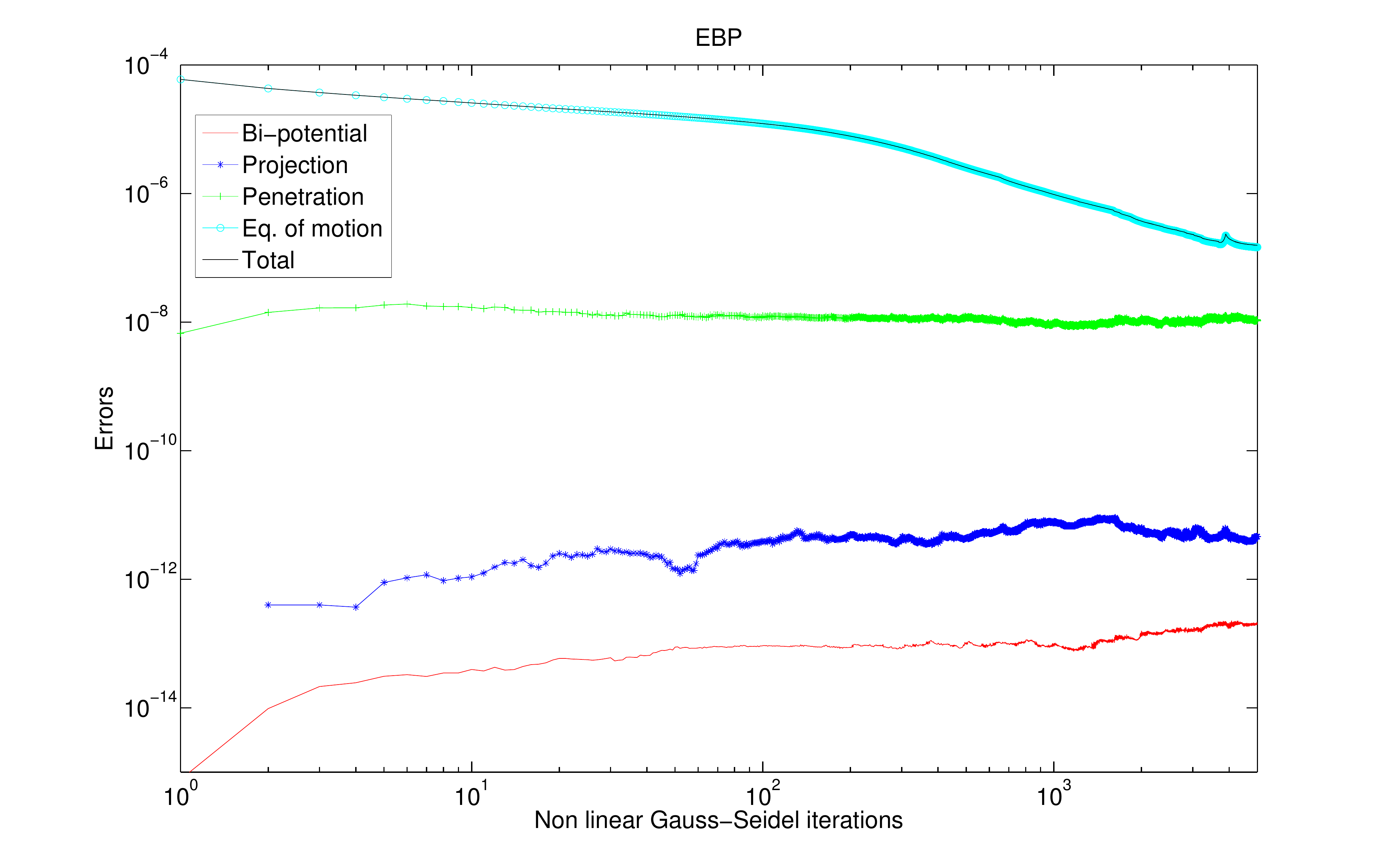}
\caption{Example 3 -- Convergence of the non-linear Gauss-Seidel iterations for the Newton and bi-potential method  (1000th time step). The two last curves collapse.}
\label{conv3-3}
\end{center}
\end{figure}

\begin{figure}[htbp]
\begin{center}
\includegraphics[width=12.0cm]{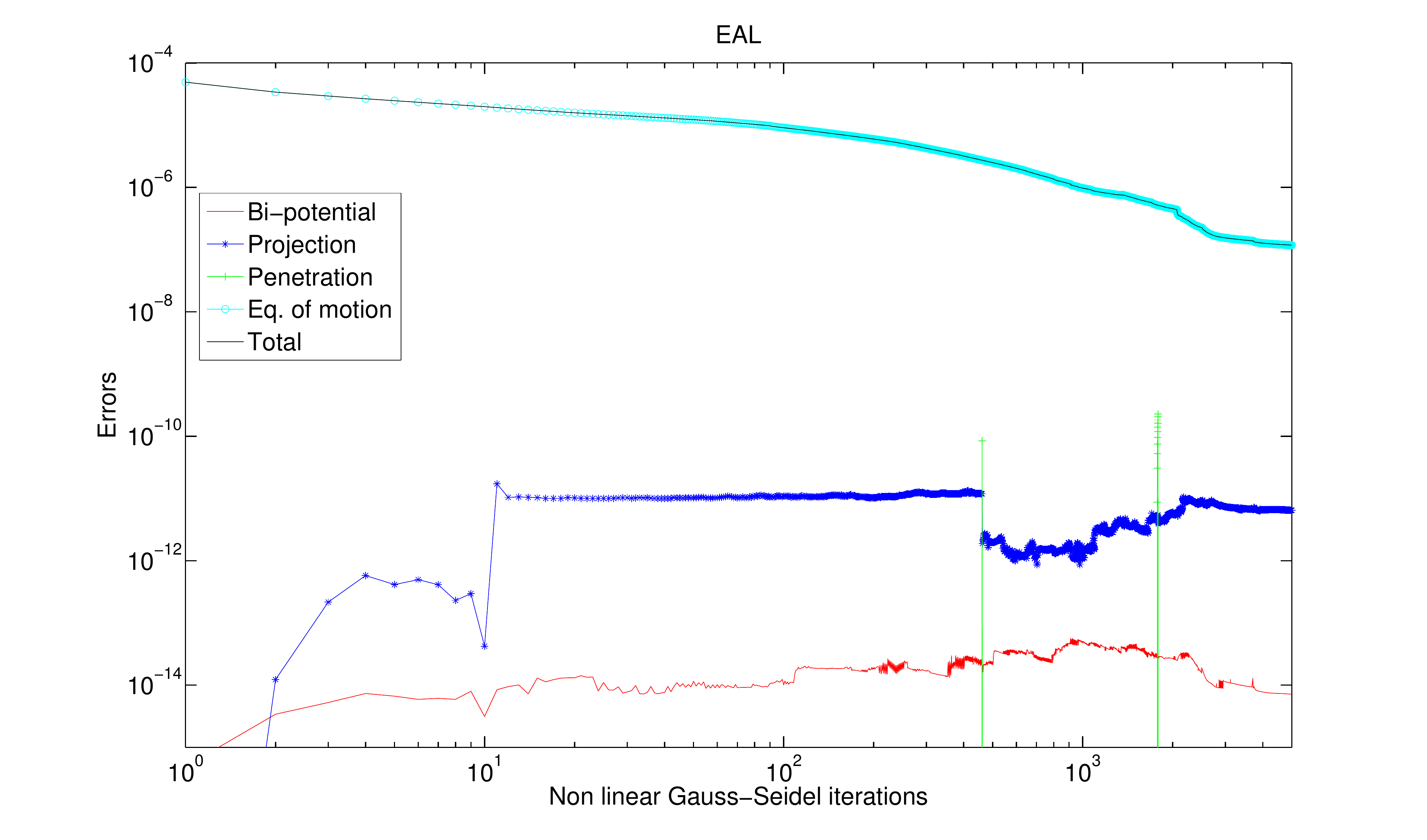}
\caption{Example 3 -- Convergence of the non-linear Gauss-Seidel iterations for the Newton and Augmented Lagragian method  (1000th time step). The two last curves collapse.}
\label{conv3-4}
\end{center}
\end{figure}

In this example, the difference between methods SBP and SAL on the one hand, and the method EBP and EAL on the other hand is larger (see table \ref{Tab:tab3-1}). 
We can notice that the SAL method is 10.95\% faster than the SBP method, and the EBP is 21.83\% faster than the SBP method. Here, the EAL method is no longer the faster one, but the penetration is very small. Again, for the two first methods the global error is essentially due to the penetrations whereas
the two last methods, the error is essentially due to the failure to follow precisely the equations of motion.

\subsection{Discussion on the descent parameter $\rho$}

We consider again the third example solved by the Newton and bi-potential method ($\varepsilon_{tot}=10^{-8}$, maximal number of iterations of Newton method equal to 100, $\varepsilon_{Newt}=10^{-5}$, $500^{th}$ time step).
Here, we take $\bar \rho=\frac{m_im_j}{m_i+m_j}\frac{1}{\Delta t}$, and we consider $\rho=\alpha\bar \rho$, for various values of $\alpha$.

\begin{table}[htbp]
\begin{center}
\begin{tabular}{|l|c|c|c|}
\hline
$\alpha$ & Number of NLGS& Maximal & Total CPU  \\
 & iterations & penetration &   time (s)  \\
 & (last time step)&  (last time step)& \\
\hline
5         &652 &$0.110\cdot 10^{-6}$ & 65.08\\
2          &414&$0.177\cdot 10^{-6}$ & 55.86\\
1         &750 &$0.149\cdot 10^{-6}$ & 53.95\\
$\frac12$ &812 &$0.634\cdot 10^{-6}$ & 78.43\\
$\frac15$ &667 &$0.219\cdot 10^{-5}$ & 176.14\\
\hline
\end{tabular}
\caption{Comparison of the results obtained for various values of $\rho=\alpha \bar \rho$ on the third example (after the 500th iteration, $N_{max}=5000$ iterations)\label{Tab:tab4-1}}
\end{center}
\end{table}%

These results show one of the main advantage of the EBP method. Indeed, one can notice that in table \ref{Tab:tab4-1}, the CPU time and the quality of the solution
are very similar if $\alpha$ is equal to 1 or 2. Even if $\alpha$ is equal to five, the convergence is not to damaged. In that case, one remain the the SBP and the SAL are no longer convergent.
If the parameter $\alpha$ is small, the method converges but the convergence rate is very small. One can notice that the EAL method is much more sensitive about
the parameter $\alpha$, essentially in the convergence of the Newton method.

\section{Conclusion}

The results presented show that, using an appropriate Newton method, it is possible to improve the computational time over that 20\%
compared to the standard methods.
Moreover, one principal drawback of that type of methods, that is the dependance of the results on the parameter $\rho$ does not exist
anymore.
\bigskip

In the future, this method will be extended to the case of a contact law with adhesion. This improvement will be realized in a near future.
\bigskip

\noindent
{\bf Acknowledgment:} This work was supported by the french CNRS grant "PEPS Opticontact", and by the project "TRIBAL" supported
by the Picardy region and the European funds FEDER.


\end{document}

%% file: signorini.pdf_t
\begin{picture}(0,0)%
\includegraphics{signorini.pdf}%
\end{picture}%
\setlength{\unitlength}{4144sp}%
\begingroup\makeatletter\ifx\SetFigFont\undefined%
\gdef\SetFigFont#1#2#3#4#5{%
  \reset@font\fontsize{#1}{#2pt}%
  \fontfamily{#3}\fontseries{#4}\fontshape{#5}%
  \selectfont}%
\fi\endgroup%
\begin{picture}(2631,1725)(259,-1063)
\put(2296,-736){\makebox(0,0)[lb]{\smash{{\SetFigFont{12}{14.4}{\rmdefault}{\mddefault}{\updefault}{\color[rgb]{0,0,0}$u_n$}%
}}}}
\put(991,-61){\makebox(0,0)[lb]{\smash{{\SetFigFont{12}{14.4}{\rmdefault}{\mddefault}{\updefault}{\color[rgb]{0,0,0}Contact}%
}}}}
\put(1891,-331){\makebox(0,0)[lb]{\smash{{\SetFigFont{12}{14.4}{\rmdefault}{\mddefault}{\updefault}{\color[rgb]{0,0,0}No contact}%
}}}}
\put(541,479){\makebox(0,0)[lb]{\smash{{\SetFigFont{12}{14.4}{\rmdefault}{\mddefault}{\updefault}{\color[rgb]{0,0,0}$r_n$}%
}}}}
\end{picture}%

%% file: coulomb1.pdf_t
\begin{picture}(0,0)%
\includegraphics{coulomb1.pdf}%
\end{picture}%
\setlength{\unitlength}{4144sp}%
\begingroup\makeatletter\ifx\SetFigFont\undefined%
\gdef\SetFigFont#1#2#3#4#5{%
  \reset@font\fontsize{#1}{#2pt}%
  \fontfamily{#3}\fontseries{#4}\fontshape{#5}%
  \selectfont}%
\fi\endgroup%
\begin{picture}(3017,2040)(328,-793)
\put(2746,1019){\makebox(0,0)[lb]{\smash{{\SetFigFont{12}{14.4}{\rmdefault}{\mddefault}{\updefault}{\color[rgb]{0,0,0}Sliding}%
}}}}
\put(1396,1064){\makebox(0,0)[lb]{\smash{{\SetFigFont{12}{14.4}{\rmdefault}{\mddefault}{\updefault}{\color[rgb]{0,0,0}$\st_t$}%
}}}}
\put(3016,119){\makebox(0,0)[lb]{\smash{{\SetFigFont{12}{14.4}{\rmdefault}{\mddefault}{\updefault}{\color[rgb]{0,0,0}${\bf u}_t$}%
}}}}
\put(991,524){\makebox(0,0)[lb]{\smash{{\SetFigFont{12}{14.4}{\rmdefault}{\mddefault}{\updefault}{\color[rgb]{0,0,0}$\mu r_n$}%
}}}}
\put(1441,-556){\makebox(0,0)[lb]{\smash{{\SetFigFont{12}{14.4}{\rmdefault}{\mddefault}{\updefault}{\color[rgb]{0,0,0}-$\mu r_n$}%
}}}}
\put(2161,-241){\makebox(0,0)[lb]{\smash{{\SetFigFont{12}{14.4}{\rmdefault}{\mddefault}{\updefault}{\color[rgb]{0,0,0}Sticking}%
}}}}
\end{picture}%

%% file: article.bbl
\begin{thebibliography}{99}
\bibitem{AC91}
P. Alart and A. Curnier.
A mixed formulation for frictional contact problems prone to
Newton like solution method. Computer Methods in Applied Mechanics and Engineering,
92(3):353–375, 1991.

\bibitem{dSF91}
G. De Saxc\'e, Z.-Q. Feng.
New inequality and functional for contact with friction: the implicit standard material approach.
Mechanics of Structures and Machines, 19:301-325, 1991.

\bibitem{FJCM05}
Z.-Q. Feng, P. Joli, J.-M. Cros, B. Magnain,
The bi-potential method applied to the modeling of dynamics problems with friction.
Math. Comput., 36:375-383 (2005).

\bibitem{F00}
J. Fortin,
Simulation numérique de la dynamique des systèmes multicorps
appliquée aux milieux granulaires.
PhD Thesis (en french), University of Lille 1, 2000.

\bibitem{FHdS02}
J. Fortin, M. Hjiaj, G. de Saxc\'e,
{\it An improved discrete element method based on a variational formulation of the contact law},
Computers and Geotechnics, 29 (2002), 609-640.

\bibitem{Fortin2005} 
J. Fortin, O. Millet, G. de Saxc\'e, 
{\it Numerical Simulation of Granular Materials by an improved Discrete Element Method},
 Int. J. for Num. Methods in Engineering, no. 62, pp. 639-663, (2005).

\bibitem{JM92}
M. Jean and J. J. Moreau.
Unilaterality and dry friction in the dynamics of rigid bodies
collection. In A. Curnier, editor, Contact Mechanics International Symposium, pages 31–
48. Presses Polytechniques et Universitaires Romanes, 1992.

\bibitem{J99}
M. Jean. 
The non smooth contact dynamics method. Compt. Methods Appl. Math. Engrg.,
177:235–257, 1999.

\bibitem{JF08}
P. Joli, Z.-Q. Feng.
Uzawa and Newton algorithms to solve frictional contact problems within the bi-potential framework.
Int. J. for Num. Methods in Engineering, no. 73, pp. 317-330, (2008).

\bibitem{JAJ98}
F. Jourdan, P. Alart, M. Jean.
A Gauss-Seidel like algorithm to solve frictional contact problems.
Computer Methods in Applied Mechanics an Engineering,
155:31-47 (1998).

\bibitem{M88}
J. J. Moreau.
Unilateral contact and dry friction in finite freedom dynamics. In J.J.Moreau
and eds. P.-D. Panagiotopoulos, editors, Non Smooth Mechanics and Applications, CISM
Courses and Lectures, volume 302 (Springer-Verlag,Wien, New York), pages 1–82, 1988.

\bibitem{M94}
J. J. Moreau. Some numerical methods in multibody dynamics: application to granular
materials. Eur. J. Mech. A Solids, 13(4-suppl.):93–114, 1994.

\bibitem{M99} 
J. J. Moreau. 
Numerical aspect of sweeping process. Compt. Methods Appl. Math. Engrg.,
177:329–349, 1999.

\bibitem{RA04}
M. Renouf and P. Alart. 
Conjugate gradient type algorithms for frictional multicontact
problems: applications to granular materials. Comput. Methods Appl. Mech. Engrg.,
194(18-20):2019–2041, 2004.

\bibitem{RA04b}
M. Renouf and P. Alart. 
Gradient type algorithms for 2d/3d frictionless/frictional multicontact
problems. In K. Majava P. Neittaanmaki, T. Rossi, O.Pironneau (eds.) R. Owen,
and M. Mikkola (assoc. eds.), editors, ECCOMAS 2004, Jyvaskyla, 24–28 July 2004.

\bibitem{RDA04}
M. Renouf, F. Dubois, and P. Alart. 
A parallel version of the Non Smooth Contact Dynamics
algorithm applied to the simulation of granular media. J. Comput. Appl. Math.,
168:375–38, 2004.

\bibitem{RAD05}
M. Renouf, V. Acary, and G. Dumont. 
Comparison of algorithms for collisions, contact
and friction in view of real-time applications. In Multibody Dynamics 2005 proceedings,
International Conference on Advances in Computational Multibody, Madrid, 21-24 June
2005.

\bibitem{RBDA05}
M. Renouf, D. Bonamy, F. Dubois, and P. Alart. Numerical simulation of two-dimensional
steady granular flows in rotating drum: On surface flow rheology. Phys. Fluids,
17:103303, 2005.

\bibitem{R}
 R.T. Rockafellar. Convex Analysis, Vol. 28, Princeton Math. Series, Princeton Univ. Press, 1970.

\bibitem{R91}
N. A. Ramaniraka, Thermom\'emcanique des contacts entre deux solides d\'eformables, EPFL PhD no. 1651 (1991),
Ecole Polytechnique F\'ed\'erale de Lausanne.

\bibitem{S06}
I. Sanni.
Modélisation et simulation bi et tri-dimensionnelles de la dynamique unilatérale des systèmes multi-corps de grandes tailles: application aux milieux granulaires. 
PhD Thesis (en french), University of Picardy, 2006.



\end{thebibliography}
